\documentclass[11pt]{article}
\usepackage{amsmath,amssymb}
\usepackage{pstricks}
\usepackage{pst-node}
\usepackage[dvips]{graphicx}
\usepackage{pst-poly}
\usepackage{multido}
\usepackage{color} 

\newtheorem{result}{Theorem}
\newtheorem{deduce}{Corollary}

\newtheorem{define}{Definition}
\newtheorem{support}{Lemma}

\newtheorem{theorem}{Theorem}

\newcommand{\qed}{%
\ifmmode 
\else \leavevmode\unskip\penalty9999 \hbox{}\nobreak\hfill \fi
\quad\hbox{\qedsymbol}}
\newcommand{\openbox}{\leavevmode \hbox to.77778em{%
\hfil\vrule
\vbox to.675em{\hrule width.6em\vfil\hrule}%
\vrule\hfil}}
\newcommand{\qedsymbol}{\openbox}

\newcommand{\showgrid}{}
\newcommand{\gridon}{\renewcommand{\showgrid}{\psset{subgriddiv=1,griddots=10,gridlabels=6pt}\psgrid}}
\gridon

\begin{document}
\begin{center}
{\bf\LARGE  Revisiting the Hamiltonian Theme in the Square of a Block: The General Case }
\end{center}

\vskip8pt

\centerline{\large   Herbert Fleischner$^a$ \  and \ Gek L. Chia$^{b,c}$}

\begin{center}
\itshape\small	
$^a$Institute of Logic and Computation, Algorithms and Complexity Group, \\ Technical University of Vienna, Austria  \\
$^b$Department of Mathematical and Actuarial Sciences, \\Lee Kong Chian Faculty of Engineering and Science,\\ Universiti Tunku Abdul Rahman, Sungai Long Campus,\\  Jalan Sungai Long, Bandar Sungai Long, \\Cheras 43000 Kajang,Selangor Malaysia \\
$^c$Institute of Mathematical Sciences, University of
Malaya, \\ 50603 Kuala Lumpur,  Malaysia  \\
\end{center}

\vspace{3mm}

\begin{abstract}
\vspace{3mm}  This is the second part of joint research in which we show that every $2$-connected graph $G$ has the ${\cal F}_4$ property. That is, given distinct $x_i\in V(G)$, $1\leq i\leq 4$, there is an $x_1x_2$-hamiltonian path in $G^2$ containing different edges $x_3y_3, x_4y_4\in E(G)$ for some $y_3,y_4\in V(G)$. However, it was shown already in \cite[Theorem 2]{cf1:refer} that 2-connected DT-graphs have the ${\cal F}_4$ property; based on this result we generalize it to arbitrary $2$-connected graphs. We also show that these results are best possible. 

 \vspace{5mm}
\end{abstract} \vspace{3mm}

\section{Introduction}

This is the second part of joint research in which we establish the most general result for the square of a block (i.e., a $2$-connected graph) to be hamiltonian connected. In the first part this was achieved in \cite[Theorem 2]{cf1:refer} for the case of DT-graphs (i.e., graphs in which every edge is incident to a vertex of degree two). In the past, the approach to deal with $2$-connected DT-graphs first and then generalize the corresponding results to blocks in general, was a logical consequence of the proof methods developed in \cite{fle1:refer}--\cite{fh:refer}, say. However, since the 1990's shorter proofs of what has become known as Fleischner's Theorem, were developed first by \v R{\'\i}ha in \cite{rih:refer}  and later by Georgakopoulos in \cite{geo:refer}. A short proof of an even stronger version of that theorem was proved by M\"uttel and Rautenbach in \cite{mr:refer}. Unfortunately, the methods developed for these shorter proofs do not seem to suffice to prove the main result of this paper (Theorem 4). This is why we had to resort to the concept of EPS-graphs (see, e.g., \cite{fle1:refer}). 

All concepts not defined in this paper, can be found in the cited literature; in cases where contradictions regarding terminology may arise, we prefer the definitions as given in the papers by Fleischner. We also included some additional references to give the interested reader a better insight regarding past developments of the topic. However, to make it easier to read this paper we repeat some definitions. In particular, by a  {\em $uv$-path} we mean a path from $u$ to $v$. If a $uv$-path is hamiltonian, we call it a {\em $uv$-hamiltonian path}. Also, we understand an {\em eulerian} graph to be a not necessarily connected graph all of whose vertices have even degree. Moreover, we let $\delta u=u$ if $d(u)=1$, and $\delta u=\emptyset$, otherwise.

Next, we repeat some results quoted or proved in \cite{cf1:refer}, using  the same numbering as in \cite{cf1:refer}. Theorems proved in the 1970's and quoted already in \cite{cf1:refer} are numbered by upper-case letters using the same letters as in \cite{cf1:refer}.
\vspace{2mm}

\begin{define}
Let $G$ be a graph and let $A = \{x_1, x_2, \ldots, x_k\}$ be a set of $k \ (\geq 3)$ distinct vertices  in $G$.  An   $x_1 x_2$-hamiltonian path in $G^2$ which contains $k-2$ distinct edges $x_iy_i \in E(G)$, $i = 3, \ldots, k$  is said to be ${\cal F} _k$. Hence we speak of an ${\cal F}_k$ $x_1x_2$-hamiltonian path in $G^2$. If $x_i$ is adjacent to $x_j$, we insist that $x_iy_i$ and $x_jy_j$ are distinct edges.  A graph $G$ is said to have the ${\cal F} _k$ property if for any set $A =\{x_1,x_2 \ldots, x_k\} \subseteq V(G)$, there is an  ${\cal F} _k$   $x_1x_2$-hamiltonian   path in $G^2$.
\end{define}

\vspace{2mm}
By an {\it EPS-graph}, {\it JEPS-graph} of $G$, denoted $S=E\cup P$, $S=J\cup E\cup P$ respectively, we mean a spannning connected subgraph $S$ of $G$ which is the edge-disjoint union of an eulerian graph $E$ (which may be disconnected) and a linear forest $P$, respectively a linear forest $P$ together with an open trail $J$.

\vspace{2mm}

\begin{support}   \label{flelemma}
([3, Lemma 1]) Suppose $G$ is a block chain with a cutvertex,  $v$ and $w$ are vertices in different endblocks of $G$ and are not cutvertices. Then

\vspace{1mm}
(i)    there exists an $EPS$-graph $E \cup P \subseteq G$ such that  $d_P(v), \ d_P(w) \leq 1$. If the endblock which contains $v$ is $2$-connected, then we have $d_P(v) =0$ and $d_P(w) \leq 1$; and

\vspace{1mm}
(ii)  there exists a $JEPS$-graph $ J \cup E \cup P \subseteq G$ such that $d_P(v) =0 = d_P(w)$.  Moreover,  $v, w$ are the only odd vertices of $J$. Also, we have  $d_P(c) = 2$ for at most one cutvertex $c$ of $G$ (and hence $d_P(c') \leq 1$ for all other cutvertices $c'$ of $G$).
\end{support}

\vspace{2mm}

\begin{theorem} \label{2edge} 
	[3, Theorem 1])
Suppose $G$ is a $2$-connected graph and $v, w$ are two distinct vertices in $G$.
Then either

(i) there exists an $EPS$-graph $S=E \cup P \subseteq G$ with  $d_P(v)=0= d_P(w)$;

or

(ii) there exists a $JEPS$-graph $S= J\cup E \cup P \subseteq G$ with $v,w$ being the only odd vertices of $J$, and $d_P(v)=0= d_P(w)$.
\end{theorem}

\vspace{3mm}





\vspace{1mm}
By a $[v; w_1, \ldots, w_n]$-$EPS$-graph of $G$, we mean an $EPS$-graph $S = E \cup P$ of $G$ such that $d_P(v)=0$ and $d_P(w_i) \leq 1$ for every $i=1, \ldots, n$.

\vspace{1mm}
\begin{theorem} 
(\cite[Theorem 3]{fh:refer}) \label{thm3fh} 
Let $G$ be a $2$-connected graph and let $v, w_1, w_2, w_3$ be four distinct vertices of $G$. Suppose $K$ is a cycle in $G$ such that $\{v, w_1, w_2, w_3\} \subseteq K$. Then $G$ has a $[v; w_1, w_2, w_3]$-$EPS$-graph $S = E \cup P$ such that $K \subseteq E$.
\end{theorem}

\vspace{1mm}
Suppose $G$ is a $2$-connected graph and $v, w_1, w_2$ are distinct vertices in $G$. A cycle $K$ in $G$ is a {\em $[v; w_1, w_2]$-maximal cycle} in $G$ if $\{v, w_1\} \subseteq V(K)$, and  $w_2 \in V(K)$ unless $G$ has no cycle containing all of $\{v, w_1, w_2\}$.

\begin{theorem} 
(\cite[Theorem 2]{fh:refer}) \label{thm2fh}
Let $G$ be a $2$-connected graph and let $v, w_1, w_2$ be three distinct vertices of $G$. Suppose $K$ is a $[v; w_1, w_2]$-maximal cycle in $G$. Then $G$ has a $[v; w_1, w_2]$-$EPS$ graph $S = E \cup P$ such that $K \subseteq E$.
\end{theorem}

\begin{theorem} 
	  (\cite[Theorem 2]{fle1:refer})   \label{thm1f}
Let $G$ be a $2$-connected graph and let $v, w$ be two distinct vertices of $G$. Let $K$ be a cycle through $v, w$. Then $G$ has a $[v; w]$-$EPS$-graph $S= E \cup P$ with $K \subseteq E$.
\end{theorem}

\begin{theorem} 
	 (\cite[Theorem 3]{fle3:refer}). Suppose $v$ and $w$ are two arbitrarily chosen vertices of a $2$-connected graph $G$.
Then $G^2$ contains a hamiltonian cycle $C$ such that the edges of $C$ incident to $v$ are in $G$ and at least one of
the edges of $C$ incident to $w$ is in $G$. Further, if $v$ and $w$ are adjacent in $G$, then these are three
different edges.  \label{fletheorem3}
\end{theorem}

\vspace{1mm}
A hamiltonian cycle in $G^2$ satisfying the conclusion of Theorem \ref{fletheorem3} is also called a $[v;w]$-hamiltonian cycle.  More generally, a hamiltonian cycle $C$ in $G^2$ which contains two edges of $G$ incident to $v$, and at least one  edge $G$ incident to each $w_i$, $i=1, \ldots, k$,  is called a
$[v; w_1, \ldots, w_k]$-hamiltonian cycle, provided the edges in question are all different.

\begin{theorem} 
	 (\cite[Theorem 4]{fle3:refer}). Let $G$ be a  $2$-connected graph. Then the following hold.

(i) $G$ has the ${\cal F} _3$ property.

(ii) For a given $q \in \{x, y\}$, $G^2$ has an $xy$-hamiltonian path containing an edge of $G$ incident to $q$.
     \label{fletheorem4}
\end{theorem}

\vspace{1mm}
By applying Theorems \ref{fletheorem3} and \ref{fletheorem4} to each block of a block chain $B$, we have the following.

\begin{deduce} \label{flecor} 
Suppose $B$ is a  non-trivial  block chain with $|V(B)| \geq 3$ and  $v$ and $w $ are vertices in different  endblocks of $G$.  Assume further that $v, w $ are not cutvertices of $B$. Then

(i) $B^2$ has a hamiltonian cycle which contains an edge of $B$ incident to $v$ and an edge of $B$ incident to  $w$. In the case that the endblock which contains $v$ is $2$-connected, then  $B^2$ has a hamiltonian cycle which contains two edges of $B$ incident to $v$ and an edge of $B$ incident to $w$.  Also,

(ii) $B^2$ has  a $vw$-hamiltonian path containing an edge of $B$ incident to $v$ and an edge of $B$ incident to $w$.
\end{deduce}


 Recall that a graph is called a {\em $DT$-graph} if every edge is incident to a $2$-valent vertex. If $G$ is a graph, we let $V_2(G)$ denote the set of all vertices of degree $2$ in $G$.

The main result of \cite{cf1:refer} is the following result which is the larger part of the proof of Theorem 4 below.

\vspace{1mm}
\begin{result}   
	\label{dt}
Every $2$-connected $DT$-graph  has the ${\cal F} _4$ property.

\end{result}

\vspace{3mm}
In proving Theorem 1 we made use of the following Lemma which plays a role also in this paper.
\vspace{2mm}

\begin{support}  \label{lemma4cd} 
Let $G$ be a $2$-connected DT-graph and let  $G^+ =G\cup\{x_1y, x_2y, y\}$, $y\not\in V(G)$ (see \cite{cf1:refer}), with $N(x_3) \not \subseteq V_2(G)$ and $N(x_4) \not \subseteq V_2(G)$. Suppose $N(x_i) \subseteq V_2(G)$ for some $i \in \{ 1,2\}$.     Assume further that every proper $2$-connected subgraph of $G$ has the ${\cal F} _4$ property. Then $(G^+)^2$ has a hamiltonian cycle containing the edges $x_1y, x_2y, x_3z_3, x_4z_4$ where $x_3z_3,x_4z_4$ are different edges of $G$. 
\end{support}

Note that in the ensuing discussion and proofs we make use of the fact that in $DT$-graphs $G$, the existence of an EPS-graph of $G$ yields a hamiltonian cycle of $G^2$. In order to keep the paper as short as possible the reader is referred to the constructions expounded in \cite {fle1:refer}.

\vspace{2mm}
However, before dealing with the main result, Theorem~4 in section~3, we need to prove several preliminary results.

\section{Beyond ${\cal F}_3$}

We now proceed to prove some results needed to shorten the proof of Theorem 4.

\begin{support}  \label{lemma2connected} 
Let $G$ be a $2$-connected $DT$-graph with at least four vertices, and let $v, w_1, w_2$ be three distinct vertices in $G$ with $N(v) \subseteq V_2(G)$ and $N(w_1) \subseteq V_2(G)$. Then $G^2$ has a $[v; w_1,w_2]$-hamiltonian cycle.  
\end{support}

\vspace{2mm}
\noindent
{\bf Proof:} Since $G$ is a $2$-connected graph, $G$ has a  cycle $K$  containing $v, w_1$. Suppose $K$ has been chosen such that it is $[v; w_1, w_2]$-maximal.  Then $G$ has a $[v; w_1, w_2]$-$EPS$-graph $S= E \cup P$ with $K \subseteq E$ by Theorem \ref{thm2fh}.  

\vspace{2mm}
If $N(w_2) \subseteq V_2(G)$, then it is straightforward to see that $S^2$ yields a hamiltonian cycle having the required properties. Note that the case $N(w_2) = \{v,w_1\}$ yields $vw_1\notin E(G)$ in this case  (since $|V(G)|\geq4$) and $K$ contains a path $v_0vw_2w_1w_0$ or $G=K=C_4$ ($v_0\in N(v), w_0 \in N(w_1)$) all of whose vertices are $2$-valent in $G$ and thus the four edges of that path are contained in some hamiltonian cycle of $S^2$.  Hence $N(w_2)\nsubseteq V_2(G)$ and $w_2$ is a $2$-valent vertex.

	

\vspace{2mm} Depending on the position of $w_2$ vis-a-vis $v$ and $w_1$ we now consider the following cases.

\vspace{2mm} {\bf Case (A)} $N(w_2)=\{v,w_1\}$. It is easy to see that $vw_1\notin E(G)$.

Next we need to consider two cases separately.

\vspace{2mm}	
(1) $G-w_2$ is $2$-connected. We apply Theorem A and correspondingly consider the following cases.

\vspace{2mm}
First we assume that $G-w_2$ has an $EPS$-graph $S=E\cup P$ with $d_P(v)=d_P(w_1)=0$. By the construction according to the method developed in [6] we have in $(G-w_2)^2$ a hamiltonian cycle $H$ whose edges in $v$ and in $w_1$ are in $G-w_2$. Now it is trivial to expand $H$ to a hamiltonian cycle in $G^2$ as required.
	 
\vspace{2mm}	
On the other hand, if $G-w_2$ has a $JEPS$-graph $S=J\cup E\cup P$ with $v, w_1$ being the only odd vertices of $J$ and $d_P(v)=d_P(w_1)= 0$, then $(G-w_2)^2$ has a hamiltonian path $P(v,w_1)$ starting in $v$ with an edge of $G$ and ending in $w_1$ with an edge of $G$, then $P(v,w_1)\cup\{w_1w_2,w_2v\}$ defines a hamiltonian cycle of $G^2$ as claimed by the lemma. 

\vspace{2mm}
(2) $G-w_2$ is not $2$-connected; hence it is a block chain with $v$ and $w_1$ belonging to different endblocks of $G-w_2$, and they are not cutvertices  of $G-w_2$. By Corollary 1{\it (ii)}, $(G-w_2)^2$ has a hamiltonian path $P(v,w_1)$ starting in $v$ with an edge of $G$ and ending in $w_1$ with and edge of $G$. Thus $P(v_1w_1)\cup\{w_1w_2,w_2v\}$ defines hamiltonian cycle of $G^2$ as claimed by the lemma and thus finishes Case (A). 

Because of the cases already treated it follows that there is $t\in V(G)$ satisfying 

\vspace{2mm}
{\bf Case (B)} $t\in N(w_2)-V_2(G)$. We assume additionally $|N(w_2)\cap\{v,w_1\}|=1$.

(i) $t\in\{v,w_1\}$. Let $t'=N(w_2)-t$. Hence $vw_1\notin E(G)$ and $t'\notin\{v,w_1\}$. Moreover, $w_2\in V(K)$ and $t'\in V_2(G)$; otherwise we could treat $t'$ like $t$ in (ii) below. In this case we can write $$K=v,\dots ,t',w_2,w_1,w'_1,\dots, v'',v\quad\quad\text{if } t=w_1$$
or
$$K= v,w_2,t',\dots,w_1,w'_1,\dots,v'',v\quad\text{if } t=v.$$
In any case, a $[v;w_1,w_2]$-$EPS$-graph $S=E\cup E$ with $K\subseteq E$ exists by Theorem~C and yields in $S^2$ a hamiltonian cycle of $G^2$ as required. 

(ii) $t\not\in\{v,w_1\}$. Since $\{v,w_1,w_2\}\subset V(K)$ we also have $t\in V(K)$, and by Theorem~C, a $[v; t, w_1]$-$EPS$-graph $S=E\cup P$ with $K\subseteq E$ exists. Also in this case, $S^2$ has a hamiltonian cycle as claimed by the lemma (in particular, it contains $tw_2$).

We are thus led to the following case.

\vspace{2mm}
{\bf Case (C)} $t\in N(w_2)-V_2(G)$ and $N(w_2)\cap\{v,w_1\}=\emptyset$. 

\vspace{2mm}
Further we assume that $w_2$ is not contained in the cycle $K$; otherwise, for $t$ as above, $K$ contains $v, w_1, w_2, t$, and  $G$ has a $[v; w_1, w_2, t]$-$EPS$-graph $S= E \cup P$ with $K \subseteq E$ by Theorem \ref{thm3fh}. Again, $S^2$ yields a hamiltonian cycle  with the required properties.

\vspace{2mm} Partition $K$ into two $vw_1$-paths, $K=P_1(v,w_1) \cup P_2(v, w_1)$.  Since $G$ is $2$-connected, there exists a $w_2u_1$-path $P(w_2, u_1)$ and a $w_2u_2$-path $P(w_2, u_2)$ in $G$ which are internally disjoint, with  $u_1, u_2 \in V(K)$ and such that $(P(w_2, u_i) -u_i) \cap K = \emptyset$, $i =1, 2$.

\vspace{2mm}
Suppose $u_1, u_2 \in P_j(v, w_1)$ for some $j \in \{1, 2\}$. Then there is a cycle $K^*$ in $G$ containing the vertices $v, w_1, w_2, t$ which contradicts the choice of $K$. 

\vspace{2mm} Hence we assume that $u_i \in P_i(v, w_1)$, $i=1, 2$;  it is an internal vertex of $P_i(v,w_1)$. 

\vspace{2mm} Now consider \[ \min  _{K \supset \{v, w_1\}} \ \min _{u_1, u_2 \in V(K)} \  \{ \ l(P(w_2, u_1)) + l(P(w_2, u_2)) \  \} \, ;         \]
 fix a cycle $K$ and $u_1, u_2 \in V(K)$ together with $P(w_2, u_1), P(w_2, u_2)$ which satisfy this minimality condition.

 \vspace{2mm} Set $ P(w_2) = P(w_2, u_1)  \cup P(w_2, u_2)$ and let $G_2 \subset G$ be induced by $ V(P(w_2))$ and by all vertices $y$ lying on a path $P_y$ with endvertices $v_y, w_y \in V( P(w_2))$ such that $\{v_y,w_y\}\not=\{u_1,u_2\}$ and satisfying $(V(P_y) - \{v_y, w_y\}) \cap V(K) = \emptyset$. $G_2$ is uniquely determined and it is a (trivial or non-trivial) block chain with $u_1, u_2$ belonging to endblocks of $G_2$; they are not cutvertices of $G_2$.

 \vspace{2mm} Likewise, define $G_K$ as induced by all vertices $z$ lying on a path $P_z$ with endvertices $v_z, w_z \in V(K)$ and satisfying $(V(P_z) - \{v_z, w_z\}) \cap V(G_2) = \emptyset$. $G_K$ is $2$-connected because of $K\subset G_K$.
 
 \vspace{2mm} Observe, that the minimality condition guarantees that there is no path $P(x,y)$ with $x\in V(G_K)-\{u_1,u_2\}$ and $y\in V(G_2)-\{u_1,u_2\}$. Now it is straightforward to see that $G= G_K \cup G_2$,    $G_K \cap G_2= \{u_1, u_2\}$ because of the minimality condition.



\vspace{2mm}
Note that  the above arguments  apply to arbitrary $2$-connected graphs. In what follows we restrict ourselves to $DT$-graphs.

\vspace{2mm}
 Also, from the choice of $K$ it follows that $\{u_1, u_2\} \cap \{v, w_1\} = \emptyset$.  However\linebreak  $K \supset \{v, w_1, u_1, u_2\}$ which is a set of four distinct vertices on $K$.  Hence $G_K$ has a $[v; w_1, u_1, u_2]$-$EPS$ graph $S_K = E_K \cup P_K$ with $K \subseteq E_K$ because of Theorem B.

\vspace{2mm}
Now consider the graph $G_2$.

\vspace{2mm} (a) Suppose $w_2$ is incident with a bridge of $G_2$. Then $G_2$ has  an $EPS$-graph\linebreak $S_2 = E_2 \cup P_2$ with $w_2 \not \in E_2$, $d_{P_2}(w_2)=2$ and $d_{P_2}(u_i) \leq 1$, $i =1, 2$ by Lemma \ref{flelemma}~(i). It follows that for $E = E_K \cup E_2$ and $P = P_K \cup (P_2 - w_2t)$, $t\in N(w_2)$, $S= E \cup P$ is an $EPS$-graph of $G$ with $K \subseteq E$, $d_P(w_2) =1$ and $w_2$ is a pendant vertex in $S$, $d_P(v) =0$, $d_P(w_1) \leq 1$, and $d_P(u_i) \leq 2$, $i=1, 2$. It now follows that $S^2$ yields a hamiltonian cycle $C$ in $G^2$ as required: its edges incident to $v$ are edges of $G$,  and at least one edge of $C$  incident to  $w_i$ is in $G$, $i=1, 2$.

\vspace{2mm}
(b) Suppose $w_2$ lies in a cycle of $G_2$, i.e., $w_2$  lies in a $2$-connected block $B(w_2)$ of $G_2$. Let $z_i \in V(B(w_2))$ be  such that $z_i = u_i$ if $u_i \in V(B(w_2))$, $i =1, 2$; otherwise, let $z_i$ be a cutvertex of $G_2$.

\vspace{2mm}
If $G_2$ is a non-trivial block chain we apply Corollary 1{\it (i)} to obtain a hamiltonian cycle $C_2$ of $G^2_2$. $C_2$ contains $u_1y_1, u_1v_1,  u_2v_2 \in E(C_2)\cap E(G_2)$ provided the endblock $B(u_1)$ containing $u_1$ is $2$-connected; and $\{y_1,v_1\}\subseteq N(u_1)$, $v_2\in N(u_2)$. However, if $B(u_1)$ is a bridge $u_1y_1$ then $u_1y_1\in E(C_2)$, and $u_1v_1\in E(C_2)-E(G)$. 
Moreover, in constructing $C_2$ (which results from applying Theorem~E to the $2$-connected blocks of $G_2$) we may apply Lemma~3 by induction to the block $B(w_2)$ containing also $s\in N(w_2)$, to obtain $w_2s\in E(C_2)$ as well. 

\vspace{2mm}
If however, $G_2$ is $2$-connected,  we apply induction to $G_2$ to obtain a hamiltonian cycle $C_2$ of $G_2^2$ where edges incident to $u_1$ are in $G_2$ and  so is $sw_2$ and an edge incident to $u_2$.

\vspace{2mm}
To obtain  $H_2$ missing $u_1$, we  make a `shortcut' by replacing $u_1y_1, u_1v_1$ with $y_1v_1$.

\vspace{2mm}
Now, $S_K$  yields a hamiltonian cycle $H_K \subseteq (G_K)^2$ with its  two edges in $v$ belonging to $G_K$ and in each of $w_1, u_1, u_2$, $H_K$ traverses at least one edge of $G_K$ (note that $N(w_1) \cup N(u_1) \cup N(u_2) \subset V_2(G)$). Likewise, $H_2$ contains an edge of $G_2$ incident with $w_2$, and one edge of $G_2$ incident with $u_2$. Denote $u_2v_K \in H_K \cap G_K$, $u_2v_2 \in H_2 \cap G_2$. Then $H = (H_K -u_2v_K) \cup (H_2 -u_2v_2) \cup \{v_Kv_2\}$ is a hamiltonian cycle $C$ in $G^2$ as required.    \qed

\vspace{3.88mm}
 By an {\em edge-critical block}, we mean a block which fails to be a block when any edge is deleted from it.

 \vspace{2mm}
Let $G$ be a graph and let $D(G) = \{ uv \in E(G) \ | \ d(u) >2, d(v) > 2 \}$. Note that  $G$ is a $DT$-graph if and only if $D(G) = \emptyset$.

\begin{theorem}  (\cite[Theorem 1]{fle3:refer})
Suppose $G$ is an edge-critical block which is not a $DT$-graph. Let $x, y$ be any two distinct vertices in $G$. Then $D(G)$ contains an edge $e$ such that $G-e$ has a $DT$-endblock $B$ such that $\{x, y \} \not \subset V(B)$, and if $x \in V(B)$, then $x$ is a cutvertex of $G-e$.   \label{fletheorem1}
\end{theorem}

\vspace{1mm} We shall now prove a  stronger version of  Theorem~F{\it (ii)}.

\begin{result}  
	\label{lemmaf3}
Let $G$ be a $2$-connected graph and let $x, y$ be two vertices in $G$. Then $G^2$ has an $xy$-hamiltonian path $P(x, y)$ such that

(i) $xz \in E(G) \cap E(P(x,y))$  for some $z \in V(G)$, 
and

(ii) either $yw \in  E(G) \cap E(P(x,y))$ for some $w\in V(G)$,  or  else $P(x, y)$ contains an edge $uv$ for some vertices $u, v \in N(y)$.
\end{result}

\vspace{2mm}  \noindent
{\bf Proof:}
Without loss of generality, assume that $G$ is edge-critical  since otherwise we can delete edges of $G$ until we reach an edge-critical block. We consider two cases.

\vspace{2mm} {\em Case (A)} $D(G) = \emptyset$.

\vspace{1mm}
Let $G^*$ denote the $2$-connected graph obtained from $G$ by adding a new vertex $z^*$ and joining $z^*$ to both $x$ and $y$.

\vspace{2mm}
First assume that  $x$ and $y$ are not adjacent in $G$.

\vspace{2mm}
(i) Assume that $N_G(x)\cup N_G(y) \subseteq V_2(G)$.

\vspace{1mm} Let $C^*$ denote any cycle containing $z^*$ and let $S^*=E^*\cup P^*$ be an $[x; y]$-$EPS$-graph of $G^*$ with $C^*\subseteq E^*$  by Theorem D. Let $S= S^* - z^*$.  Then $S= J \cup E \cup P$ is a $JEPS$-graph of $G$ with $P=P^*$ and the component of $S^*$ containing $C^*$ becomes the open trail $J$ from $x$ to $y$ in $S$. By following the construction of an $xy$-hamiltonian path $P(x, y)$ in $S^2$ which was used in \cite{fle1:refer},  
it is clear that $P(x, y)$ can start with an edge of $S$ incident to $x$ and ends with an edge of $S$ incident to $y$ unless $d_P(y) =1$.  If $d_P(y)=1$, we jump from a vertex $u$ preceding $y$ in $J$ to the vertex $v$ in $P_0$ adjacent to $y$, where $P_0$ is the component of $P$ containing $y$.

\vspace{2mm}
(ii) Assume that $N_G(x) \not \subseteq V_2(G)$ and $N_G(y) \subseteq V_2(G)$.

\vspace{1mm} Then at least one of the two neighbors of $x$, say $x'$ has degree greater than $2$. Let $C^*$ be a cycle containing $z^*$ and the edge $xx'$. Note that this is possible because $G$ is  $2$-connected (so that there is an $xy$-path in $G$ starting with any given edge). In this case, let $S^* = E^* \cup P^*$ be an $[x; x', y]$-$EPS$-graph of $G^*$ with $C^*\subset E^*$  by Theorem \ref{thm2fh}  because $C^*$ is $[x; x',y]$-maximal. Then proceed as in case (i) and note that $x$ is a pendant vertex in $S$. A  required hamiltonian path in $S^2$ (with  $S= E\cup P \cup J$  as in case (i)) can be constructed starting with the pendant edge incident to $x$.

\vspace{2mm}
(iii) Assume that $N_G(y)  \not \subseteq V_2(G)$ and $N_G(x) \subseteq V_2(G)$.

\vspace{1mm} This case can be treated symmetrically to case (ii), starting with an $[x;y',y]$-EPS-graph $S^*=E^*\cup P^*$ of $G^*$ and $y'\in (N_G(y)-V_2(G))\cap V(C^*)$.

 \vspace{2mm}
(iv) Assume that $N_G(x)  \not \subseteq V_2(G)$ and $N_G(y) \not \subseteq V_2(G)$.

\vspace{1mm} Proceed as in case (ii) with $C^*$ as defined there. Here we operate with an $[x; x', y']$-$EPS$-graph $S^* = E^* \cup P^*$ of $G^*$ with $C^* \subseteq E^*$, where  $y'y \in E(C^* -z^*)$,  assuming first that $x'\not= y'$ (i.e., $\ell(C^*)>4$) and applying Theorem C.   Then $d_{P^*}(y) \leq 1$ (because $d_{G^*}(y) =3$).  Again we get a required $xy$-hamiltonian path $HP$   in $G^2$. 
 \vspace{2mm}

Note that, if $y'' \in N_G(y)-y'$ and $d_{P^*}(y'') =2$, $d_{P^*}(y)=1$  then $yy''$ is an end-edge of the path in $P^*$ incident to $y$ and $y'y''\in E(HP)$.

\vspace{2mm}
Now assume that $x'=y'$, (i.e. $\ell(C^*)=4$). 
Since $G$ is $2$-connected, there is an $x'y$-path $P(x',y)$ in $G-x$  not containing $x'y$ ($x'y$ lies in a $2$-connected block of $G-x$).     Then $\{xx'\} \cup P(x', y)$ is an $xy$-path in $G $ which together with $xz^*y$ yields a cycle $C'\subset G^*-x'y$ with $\ell(C')>4$,   for which the preceding argument goes through if we operate with an $[x; x', y'']$-$EPS$-graph of $G^* - x'y$  where $y''$ is as above ($x'y$ is a chord of $C'$ in $G^*$).

\vspace{2mm}
Next we assume that  $x$ and $y$ are adjacent. In this case, we take a longest $xy$-path in $G-xy$ and  combine it with $xz^*y$ to form the cycle $C^*$; $l(C^*) \geq 5$ follows unless $N(x)\cap N(y)\neq\emptyset$ in which case $G=K_3$ since $G$ is a $DT$-graph and we are done. If $\ell(C^*)\geq 5$ we proceed as before.

\vspace{2mm} {\em Case (B)} $D(G) \neq \emptyset$.

\vspace{2mm}
   By \cite[Theorem 1]{fle2:refer},  $D(G)$ contains an edge  $e=st$ such that $G-e$ is a block chain  with  at least one of its endblocks, say $B_e$, being a $DT$-block. Without loss of generality $t\in V(B_e)$. 
\vspace{2mm}
   
Suppose $(V(B_e) - c_e) \cap \{x,y\} = \emptyset$, where  $c_e$ is the cutvertex of $G - e$ belonging to  $B_e$.  
Then we replace $B_e$ by a path  $P^*$ of length $3$ joining $t$ and $c_e$. The resulting graph $H$ is an edge-critical block and $|D(H)| < |D(G)|$. By induction $H^2$ has an $xy$-hamiltonian path with properties {\it (i)} and {\it (ii)} as stated by the theorem. Assuming that it contains as many edges of $H$ as possible,   any such $xy$-hamiltonian path in $H^2$ can  be converted into an $xy$-hamiltonian path in $G^2$ having properties {\it (i)}  and {\it (ii)} of the theorem,  by the same method  used  in \cite{fle2:refer} { as long as $c_e \notin\{x,y\}$. The same conclusion can be drawn if said hamiltonian path in $H^2$ satisfies $c_e\in\{x,y\}$. For, we may proceed as in \cite[ pp.~32-33]{fle2:refer}, cases 2 and 4: we just look at the $xy$-hamiltonian path 
$$ P(x,y)=c_e\dots u^*v^*\dots r\quad (u^*\in V(P^*), v^*\in V(H)-V(P^*),  \{r,c_e\}=\{x,y\})$$
in $H^2$ just as we would look at a hamiltonian cycle $H_1$ in $H^2$ in [7]
$$H_1=c_e\dots u^*v^*\dots c_e$$
and using a hamiltonian path in $B_e^2$ starting in $t$ and ending at $c_e$ with an edge of $B_e$.

\vspace{2mm} Hence we assume that for every $DT$-endblock $B_e$ of $G-e$ (where $e\in D(G)$), 
$$|(V(B_e) - c_e) \cap \{x,y\}| = 1$$
(note that $D(G) \neq \emptyset$ implies that $G$ has at least two $DT$-endblocks like $B_e$).  In particular, we assume $x \in V(B_e)-c_e$. 
   
\vspace{2mm} Let $B_e'$ be the other endblock of $G-e$. If $B_e'$ is a $DT$-block, then it follows from the preceding argument that $|(V(B_e') - c_e') \cap \{x,y\}| = 1$ where $c_e'$ is the cutvertex of $G-e$ belonging to  $B_e'$. 
   If $B_e'$ is not a $DT$-endblock, then $B_e'$ contains a $DT$-endblock $B_{e'}$ for some $e' \in D(G)$, and we have the same conclusion as in the preceding sentence. Thus we conclude in any case  that  $y \in V(B'_e)-c_e'$.
   
\vspace{2mm} Set $G_0 = G-e - (B_e \cup B'_e)$; $G_0$ is a (trivial or non-trivial) block chain. Possibly $G_0 = \emptyset$ in which case $c_e = c_e'$.
   
\vspace{2mm} By Theorem F{\it (ii)}, $(B_e)^2$ has an $xc_e$ hamiltonian path $P(x, e_c)$ starting with an  edge $xz_1$ of $B_e$; $(B_e')^2$ has an $c_e'y$-hamiltonian path $P(c_e', y)$ ending with an edge $z_2y$ of $B_e'$. 
 By Corollary 1 {\it (ii)}, $(G_0)^2$ has a $c_ec_e'$-hamiltonian path $P_0(c_e, c_e')$, being just a vertex if $c_e = c_e'$.  Then $$  P(x, c_e) P_0(c_e, c_e') P(c_e', y) $$ is an $xy$-hamiltonian path in $G^2$ having properties (i) and (ii) of the theorem.     \qed

\vspace{2mm}

\begin{define} 
  A graph  $G$ is said to have the  strong ${\cal F} _3$ property if for any set of  three distinct vertices $\{x_1, x_2, x_3\}$ in $G$, there is an  $x_1x_2$-hamiltonian   path in $G^2$ containing $x_3z_3, x_iz_i$ which are distinct edges of $G$ for a given  $i \in \{1, 2\}$. Such an $x_1x_2$-hamiltonian path in $G^2$ is called a  strong ${\cal F} _3$ $x_1x_2$-hamiltonian path.
\end{define}

\vspace{1mm}

\begin{result} \label{strongf3} 
Every  $2$-connected graph has the strong ${\cal F}_3$ property.
\end{result}

\vspace{2mm}  \noindent
{\bf Proof:} Let $G$ be a $2$-connected graph. Without loss of generality, assume that $G$ is an edge-critical block; otherwise we delete edges from $G$ until we reach an edge-critical block. Trivially, the theorem is true if $G$ is a triangle. Thus we assume that $|V(G)|\geq 4$.

\vspace{2mm} (I)  Assume that $G$ is a $DT$-graph.

\vspace{2mm}
 Proceeding analogously to what we did in proving  (\cite[Theorem 2]{cf1:refer}),  let  $G^+$ denote the graph obtained from $G$ by adding a new vertex $z$ and join $z$ to $x_1, x_2$. We shall show that $(G^+)^2$ has a hamiltonian cycle $C_i$ containing $zx_1, zx_2, x_iz_i, x_3z_3$ which are distinct edges of $G^+$ for a given $i \in \{1, 2\}$. Then $C_i - z = P_i (x_1, x_2)$ is a required strong ${\cal F}_3$ $x_1x_2$-hamiltonian path in $G^2$ containing the edges $x_iz_i, x_3z_3$ of $G$.  Basically, we apply the construction of a hamiltonian cycle in the square of an $EPS$-graph in a $DT$-graph (see [6] and  Observation (*) in [3]). In some of the cases, however, we shall  proceed by induction, noting that the theorem  is trivially true if it is a cycle; and sometimes we proceed by a  direct proof.

\vspace{2mm}
 Let $C^+$ be a cycle in $G^+$ containing $z, x_1, x_2, x_3$.

\vspace{2mm}
{\em Case (A):}    $N(x_j) \subseteq V_2(G)$, $j=1, 2, 3$.

\vspace{2mm}  By Theorem C, let $S= E \cup P$ be an $[x_i; x_{3-i}, x_3]$-$EPS$-graph of $G^+$ with $C^+ \subseteq E$. Hence $(G^+)^2$ has an $[x_i; x_{3-i}, x_3]$-hamiltonian cycle $C_i$ for any $i \in \{1, 2\}$ provided $\ell(C^+)>4$ (see the corresponding argument in the proof of Theorem 2). 

\vspace{2mm}
 However, if $\ell(C^+)=4$, then  $G^-=G-x_3$ is a non-trivial block chain ($x_1x_2\in E(G)$ yields $G$ being a triangle, contrary to the assumption at the beginning of the proof). 

\vspace{2mm}	
 Moreover, $x_1$ and $x_2$ are pendant vertices of $G^-$. By Corollary 1{\it (ii)}, $(G^-)^2$ has an $x_1x_2$-hamiltonian path $P^-_{1,2}$ starting with $x_1v_1\in E(G)$ and ending with $v_2x_2\in E(G)$. Thus
$$(P^-_{1,2}-x_1v_1)\cup\{x_1x_3, x_3v_1\}$$
and $$(P^-_{1,2}-v_2x_2)\cup\{x_3x_2,x_3v_2\}$$
yield the hamiltonian paths in $G^2$ as required by the theorem.	

\vspace{3mm}
{\em Case (B):}    $N(x_i) \subseteq V_2(G)$, $i=1, 2$ and $N(x_3) \not \subseteq V_2(G)$.

\vspace{2mm} Then $d_G(x_3) = 2$. Let $N(x_3) = \{u_3, v_3\}$.

\vspace{2mm} \hspace{5mm}
(a) $\{u_3, v_3\} \neq \{x_1, x_2\}$. Without loss of generality assume that $u_3 \not \in \{x_1, x_2\}$.

\vspace{2mm} Again, by Theorem C, let  $S = E \cup P$ be an $[x_i; x_{3-i},  u_3]$-$EPS$-graph of $G^+$ with $C^+ \subseteq E$. A required hamiltonian cycle $C_i$ in $(G^+)^2$ can be constructed using $S$.

\vspace{2mm} \hspace{5mm}
(b) $\{u_3, v_3\} =  \{x_1, x_2\}$.

\vspace{1mm}
Consider the graph $G' = G- x_3$.

\vspace{1mm} \hspace{8mm} (b1)
Suppose $G'$ is $2$-connected.  We apply Theorem A with $x_1, x_2$ in place of $v, w$.

 \vspace{1mm} (i) Suppose $G'$ has an $EPS$-graph $S' =  E' \cup P'$ with $d_{P'}(x_i) = 0$, $i=1, 2$. Let $H'$ be a hamiltonian cycle of  $(S')^2$: 
   the edges of $H'$ incident to  $x_i$, $i=1, 2$ are in $G'$; denote them by $e_i = x_iu_i, f_i = x_iv_i$, $i =1, 2$. Without loss of generality  the notation is chosen in such a way that $P(e_1, e_2)$ is the path in $H'$ starting in $x_1$ with $e_1$ and ending in $x_2$ with $e_2$; $P(f_2, f_1) \subset H'$ is defined analogously. Then 
$$  x_1 (P(e_1, e_2) -e_2)u_2v_2(P(f_2, f_1) - \{f_1, f_2\}) v_1x_3x_2   $$ 
is a hamiltonian path as required for $i=1$. By a symmetrical argument one obtains a  hamiltonian path ending with $f_2$, say, and containing $x_1x_3$.

 \vspace{1mm} (ii) Suppose $G'$ has a $JEPS$-graph $S' = J' \cup E' \cup P'$ with $x_1, x_2$ being the only odd vertices of $J'$ and $d_{P'}(x_i) = 0$,  $i=1, 2$. 
 $(S')^2$ contains a hamiltonian path                        $P^*$ starting with $g_1 = x_1y_1$ and ending with $g_2 = x_2y_2$, $\{ g_1, g_2\}  \subseteq E(G)$. We extend $P^*$  to a hamiltonian path $P$ as required by setting   $ P = x_1x_3y_1(P^*- g_1)$ or $P = (P^* - g_2)y_2x_3x_2$.

\vspace{1mm} \hspace{8mm} (b2)
Suppose $G'$ is not $2$-connected.  Then $G'$ is a block  chain. By Lemma 1{\it (ii)} with $x_1 =v$ and $x_2 =w$, $G'$ has a $JEPS$-graph $S' = J' \cup E' \cup P'$  with $d_{P'}(x_1) =0= d_{P'} (x_2)$; and $x_1, x_2$ are the odd vertices of $J'$.  Now proceed  as in (b1)(ii): $(S')^2$ has an $x_1x_2$-hamiltonian path $P^*$ starting and ending with edges $h_1, h_2$ of $G$; one extends $P^*$ to a corresponding hamiltonian path in $G^2$ by either traversing $x_1x_3$ first and ending with $h_2$ in $x_2$, or traversing $h_1$ first and ending in $x_2$ with $x_3x_2$.

\vspace{3mm}
{\em Case (C):}    $N(x_1) \subseteq V_2(G)$ and $N(x_2) \not \subseteq V_2(G)$.

\vspace{2mm} Then $d_G(x_2) =2$;  let $N(x_2) = \{u_2, v_2\}$. Without loss of generality assume that $u_2$ is on the cycle $C^+$.

\vspace{2mm}
(1) $N(x_3) \subseteq V_2(G)$.

\vspace{2mm} \hspace{5mm} (a) Suppose $x_3 \not \in N(x_2)$.

\vspace{1mm}
By Theorem B, let $S = E \cup P$ be an $[x_i; x_{3-i}, u_2, x_3]$-$EPS$-graph of $G^+$ with $C^+ \subseteq E$. Then a required hamiltonian cycle in $(G^+)^2$ can be constructed for each $i \in \{1, 2\}$.

\vspace{2mm} \hspace{5mm} (b) Suppose $x_3 \in N(x_2)$; that is, $x_3 =u_2$.

\vspace{1mm}  Let $G' = G^+ - x_2x_3$ which is a $DT$-graph.

\vspace{2mm}
(i) Suppose $G'$ is $2$-connected.

\vspace{1mm}
 There is a cycle $C'$ in $G'$ containing $z, x_1, x_2, v_2, x_3$: this follows from the fact that $G'$ contains in this case a path $P(x_1,v_2)$ with $x_3\in V(P(x_1,v_2))$; it cannot contain $z$ because $d_{G'}(x_2)=d_{G'}(z)=2$.  By Theorem C, let $S = E \cup P$ be an $[x_1; v_2, x_3]$-$EPS$-graph of $G'$ with $C' \subseteq E$. Hence a required hamiltonian cycle $C_i$ of $(G^+)^2$ can be constructed. A~corresponding hamiltonian path in $G^2$ starts and ends with edges of $G$.

\vspace{2mm}
(ii) Suppose $G'$ is not $2$-connected.

\vspace{1mm} Then $G'$ is a block chain with a $2$-connected  endblock $B_z$ containing $z, x_1, x_2, v_2$,  and a block chain $G_3 = G' - B_z$ containing $x_3$ (which is not a cutvertex of $G_3$ and belongs to an endblock of $G_3$). $G_3$ is a $DT$-graph unless $G_3=K_2$. Denote  $V(B_z) \cap V(G_3) = \{c\}$.

\vspace{1mm} By Lemma 1{\it (i)}, if $G_3$ has a cutvertex, then it has an $EPS$-graph $S_3 = E_3 \cup P_3$ such that $d_{P_3}(x_3) \leq 1$ and $d_{P_3}(c) \leq 1$.  Moreover,   if the endblock $B_c$ in $G_3$ containing $c$ is $2$-connected, then we may achieve $d_{P_3}(c)=0$; if $B_c$ is a bridge, then $d_{P_3}(c) =1 $ and $c$ is a pendant vertex in $S_3$. However, if $G_3$ is $2$-connected, then we apply Theorem D to obtain such $S_3$. {If $G_3=K_2$, then $S_3=\{cx_3\}$ and $E_3=\emptyset$, $P_3=\{cx_3\}$ and $d_{P_3}(c)=d_{P_3}(x_3)=1$.

\vspace{1mm}
Let $C_z$ be a cycle in $B_z$ containing $z, x_1, x_2, v_2, c$. Such $C_z$ exists because $d(x_2)=2$.

\vspace{1mm} If $c \not \in \{x_1, v_2\}$, then by Theorem B, let $S_z = E_z \cup P_z$ be an $[x_1; v_2, c, z]$-$EPS$-graph of $B_z$ with $C_z \subseteq E_z$. Set $E = E_z \cup E_3$ and $P = P_z \cup P_3$. Then we have an $EPS$-graph $S = E \cup P$ in $G$ with $C_z \subseteq E$ and $d_P(x_1) =0 = d_P(x_2) $, $d_P(v_2) \leq 1$ and $d_P(c) \leq 2$,  $d_P(x_3)\leq 1$. Thus a  hamiltonian cycle in $(G^+)^2$ can be constructed which contains edges of $G$ incident with $x_1, x_2$ together with another edge of $G$ incident to $x_3$,and also containing $zx_1, zx_2$.

\vspace{1mm} If $c = v_2$, then $v_2$ is a cutvertex of $G$:  for, $(G_3\cup\{cx_2,x_2x_3\})\cap (B_z-\{cx_2\}-z)=c$ and $(G_3\cup\{cx_2,x_2x_3\})\cup(B_z-\{cx_2\}-z)=G$. This yields a  contradiction.

\vspace{1mm} Hence we are left with the case $c = x_1$.

\vspace{1mm} Suppose $d_{G_3}(x_1) > 1$. Then $G_3$ has an $[x_1; x_3]$-$EPS$-graph $S_3 = E_3 \cup P_3$ which we combine with an $[x_1; v_2]$-$EPS$-graph $S_z = E_z\cup P_z$ of $B_z$ (see Theorem D) to obtain the  $EPS$-graph $S = E \cup P$ by putting $E = E_3 \cup E_z$ and $P = P_3 \cup P_z$.  We have $d_P(x_1) =0 = d_P(x_2) $, $d_P(x_3) \leq 1$, $d_P(v_2) \leq 1$, and since $G'$ is a $DT$-graph, $S^2$ has a hamiltonian cycle as required containing $x_1w_1, x_2v_2, x_3w_3 \in E(G)$,  and also containing $zx_1, zx_2$.

\vspace{1mm} Finally, assume $d_{G_3}(x_1) =1$; i.e., $G_3$ is a non-trivial block chain  or $G_3=K_2$.  Suppose first that $G_3\not= K_2$. 
By Corollary 1, $(G_3)^2$ has a hamiltonian cycle $H_3 \supset \{x_1y_1, x_3w_3 \} $ with $\{x_1y_1, x_3w_3 \} \subset E(G_3)$; and it has a hamiltonian path $P_{1,3}$ starting with $x_1y_1$ and ending with $x_3z_3$, say, which are edges of $G_3$. Likewise, since $G_2 = (G - x_2x_3)- G_3  = B_z-z$ is a non-trivial block chain ($x_2$ is a pendant vertex of $G_2$), $(G_2 - \delta x_1)^2$ has a hamiltonian cycle $H_2$ containing $v_2x_2$ and if $d_{G_2}(x_1) > 1$,  then $x_1s_1, x_1t_1 \in E(G_2)$ unless $G-\delta x_1=K_2$ in which case $H_2=x_2v_2$. It also has a hamiltonian path $P_{1,2}$ starting with $x_1z_1$, say, and ending with $v_2x_2$ which are edges of $G_2$. 

\vspace{1mm} Setting $H_2' = H_2$ if $\delta x_1 =x_1$ and $H_2' = (H_2 - \{x_1s_1, x_1t_1\}) \cup \{s_1t_1\}$ if $\delta x_1 = \emptyset$, we obtain hamiltonian paths  $P_1 (x_1, x_2)$,    
$P_2 (x_1, x_2)$ in $G^2$ as required and defined by
$$ E(P_1 (x_1, x_2)) =  E(P_{1,3} \cup (H_2' - x_2v_2)) \cup \{ x_3v_2\} $$
$$ E(P_2 (x_1, x_2)) =  (E(P_{1,2} \cup H_3)  - \{x_1z_1, x_1y_1\}) \cup \{ y_1z_1\}.$$

If, however, $G_3=K_2$, then $N_{G^+}(z)=N_{G^+}(x_3)=\{x_1,x_2\}$, i.e., $G^+-x_3$ is isomorphic to $G$. Since $d_G(x_2)=2$ and because of the assumption $d_{G_3}(x_1)=1$ and because $c=x_1$, it follows that $G-x_3$ is a non-trivial block chain and $x_2$ is an endvertex of $G-x_3$ and $c=x_1$ is not a cutvertex belonging to the other endblock of $G-x_3$ unless $x_1x_2\in E(G-x_3)$.  However, if $x_1x_2\in E(G-x_3)$, then we conclude that $G$ is a triangle in this exceptional case, contradicting the assumption $|V(G)|\geq 4$ at the beginning of the proof. Hence $G-x_3$ is a non-trivial block chain.
	
Now we apply Corollary 1{\it (ii)} to obtain in $(G-x_3)^2$ a hamiltonian path $P^{(3)}(x_1,x_2)$ starting with $x_1v_1\in E(G)$ and ending in $x_2$ with $v_2x_2\in E(G)$. Now, for $i=1,2$,
	$$(P^{(3)}(x_1,x_2)-x_iv_i)\cup \{x_ix_3,x_3v_i\}$$
yields a hamiltonian path in $G^2$ as required. 
	
\vspace{2mm}
(2) $N(x_3) \not \subseteq V_2(G)$.

\vspace{2mm} Then $d_G(x_3) = 2$. Let $N(x_3) = \{u_3, v_3\}$. Suppose without loss of generality that  $C^+$ is of the form $zx_1u_1 \ldots u_3x_3v_3 \ldots u_2x_2z$.

\vspace{2mm} \hspace{5mm}
(a) $\{u_3, v_3\} \neq \{x_1, x_2\}$.

\vspace{2mm} By Theorem B, let $S = E \cup P$ be an $[x_i; x_{3-i},  x_3^*, u_2]$-$EPS$-graph of $G$ with $C^+ \subseteq E$, where $x_3^* \in \{u_3, v_3\} - \{x_1, x_2\}$.  If $\ell(C^+)>5$, then it is straightforward to see that a~required hamiltonian cycle $C_i$ in $(G^+)^2$ can be constructed from $S$ for any $i \in \{1, 2\}$, independent of the size of $N(x_3) \cap (N(x_1) \cup N(x_2))$.

\vspace{2mm} 
Observe that $\ell(C^+)\geq 4$. However, $\ell(C^+)=4$ implies $N(x_3)=\{x_1,x_2\}$, contrary to the assumption $\{u_3,v_3\}\not=\{x_1,x_2\}$.
\vspace{2mm} 

To finish this case (a) we are thus left with the case $\ell(C^+)=5$ which implies $|\{u_3,v_3\}\cap \{x_1,x_2\}|=1$. More precisely, we have 
$$C^+=zx_1u_1u_2x_2z,$$
i.e., $x_3\in \{u_1,u_2\}$. 
\vspace{2mm} 

Suppose $d_G(u_2)=2$; then $G^-=G-\{u_1,u_2\}$ is a non-trivial block chain (the case $G^-=K_2$ is impossible).  
By Corollary 1{\it (ii)}, $(G^-)^2$ has a hamiltonian path $P^-$ starting with $x_1t_1\in E(G)$ and ending with $v_2x_2\in E(G)$. Now 
$$(P^- -x_1t_1)\cup\{x_1u_2, u_2u_1,u_1t_1\}$$
and $$(P^--v_2x_2)\cup\{v_2u_2, u_2u_1, u_1x_2\}$$
yield the required hamiltonian paths.
\vspace{2mm} 

Finally, suppose $d_G(u_2)>2$. Then $x_3=u_1$ since $d_G(x_3)=2$. 

\vspace{2mm} 
Now $G-x_3$ is either a non-trivial block chain or it is $2$-connected. In any case, $x_2$ and $u_2$ are not cutvertices of $G-x_3$ and they belong to the same $2$-connected block $B^*$ of $G-x_3$. If $G-x_3$ is not $2$-connected, let $c^*$ denote the cutvertex of $G-x_3$ in (the endblock) $B^*$. Set $G^*=(G-x_3)-B^*$. By Corollary 1{\it(ii)} or if $G^*=K_2$, $(G^*)^2$ has an $x_1c^*$-hamiltonian path $P^*$ starting with an edge $x_1t_1\in E(G)$, provided $G^*\not =\emptyset$; if $G^*=\emptyset$ set $P^*=\emptyset$. In any case, however, $(B^*)^2$ has by induction $c^*x_2$-hamiltonian paths, one starting in $c^*$ with $c^*t^*\in E(B^*)$, whereas the other ends in $x_2$ with $t_2x_2\in E(B^*)$, and both containing an edge  $u_2w_2\in E(B^*)$. Denote these paths by $P_1^*$ and $P_2^*$, respectively. If $G^*=\emptyset$, then set $c^*=x_1$.

\vspace{2mm} 

It follows that for both $i=1,2$,
$$P^*\cup(P^*_i-u_2w_2)\cup\{w_2u_1, u_1u_2\}$$
yield $x_1x_2$-hamiltonian paths as required. This finishes case (a).

\vspace{2mm} \hspace{5mm}
(b) $\{u_3, v_3\} =  \{x_1, x_2\}$.

\vspace{1mm}
Then $G' = G^+ -x_3$ is $2$-connected and it contains a cycle $C' \supseteq  \{z, x_1, x_2, v_2\}$.  By Theorem~D, let $S'= E' \cup P'$ be an $[x_1; v_2]$-$EPS$-graph of $G'$ with $C' \subseteq E'$. Then $d_{P'}(z)=0 = d_{P'}(x_1) = d_{P'}(x_2)$ and $d_{P'}(v_2) \leq 1$. Let $E = E'$ and $P = P' \cup \{x_ix_3\}$. Then $S= E \cup P$ is an $EPS$-graph of $G^+$ and a required hamiltonian cycle $C_i$ in $(G^+)^2$ containing $x_ix_3, x_{3-i}z_{3-i}$ ($z_{3-i} \in N(x_{3-i})$), which are edges of $G$, can be constructed for each $i \in \{1, 2\}$.

\vspace{2mm} Since the case $N(x_1) \not \subseteq V_2(G)$ and $N(x_2) \subseteq V_2(G)$ is symmetrical to {\it Case (C)} just considered, we are left with the consideration of one more large case for $DT$-graphs.

\vspace{3mm}
{\em Case (D):}    $N(x_1) \not \subseteq V_2(G)$ and $N(x_2) \not \subseteq V_2(G)$.

\vspace{2mm} Then $d_G(x_i) =2$ for $i = 1, 2$. Let $N(x_i) = \{u_i, v_i\}$, $i \in \{1, 2\}$. Suppose $C^+$ is of the form $zx_1u_1 \ldots u_3x_3v_3  \ldots u_2x_2z$ as in {\it Case (C)} (2) above.

\vspace{2mm}
(1)  Suppose $N(x_3) \subseteq V_2(G)$.

\vspace{2mm}
(a)  Suppose $x_3 \notin (N(x_1) \cup  N(x_2))\cap V(C^+)$; That is, $x_3 \notin \{u_1, u_2\}$.

\vspace{1mm} Let $S= E \cup P$ be an $[x_i; u_1, u_2, x_3]$-$EPS$-graph of $G^+$ with $C^+ \subseteq E$ which exists by Theorem~B. Then a required  hamiltonian cycle $C_i$ in $(G^+)^2$ can be constructed for each $i \in \{1, 2\}$.

\vspace{2mm}
(b) Suppose $x_3 \in N(x_1) \cap N(x_2)\cap V(C^+)$; hence $x_3=u_1=u_2$.

\vspace{2mm}
Suppose first that $d_G(x_3)>2$. Consider $G'_2=G-x_2x_3$. Note that $x_2$ is a pendant vertex in $G'_2$. Let $B_2$ be the endblock of $G'_2$ with $x_1,x_3\in V(B_2)$, and they are not cutvertices. Then $G'_2-B_2\neq\emptyset$ is a trivial or non-trivial block chain with $c_2=V(B_2)\cap V(G'_2-B_2)$ being the cutvertex of $G'_2$ in $B_2$. Using induction on $B_2$ we have a hamiltonian path $P(x_1,c_2)$ in $(B_2)^2$ starting with $x_1s_1\in E(B_2)$ and containing another edge $x_3s_3\in E(B_2)$. Moreover $(G'_2-B_2)^2$ has a hamiltonian path $P(c_2,x_2)$ starting with $x_2s_2\in E(G'_2-B_2)$ by Corollary 1{\it (ii)} or if $G'_2-B_2=K_2$. Then $P(x_1,c_2)P(c_2,x_2)$ is a hamiltonian path as required.

\vspace{2mm}
If however, $d_G(x_3)=2$, then we set $G''=G-x_3$ which is a non-trivial blockchain with pendant vertices $x_1,x_2$, otherwise $G=K_3$ and this case has been solved at the beginning of the proof. Thus $(G'')^2$ has a hamiltonian path $P(x_1,x_2)$ starting and ending with edges in $G''$, by Corollary 1{\it (ii)}. Now it is trivial to enlarge $P(x_1,x_2)$ to a hamiltonian path $P$ of $G^2$ as required by appropriately using $x_3x_i$, $i\in\{1,2\}$ as the last edge in $P$. This finishes case (b).

\vspace{2mm}
(c) Suppose  $x_3 \in (N(x_i)-N(x_{3-i}))\cap V(C^+)$, $i\in\{1,2\}$. Without loss of generality $i=2$. Hence $x_3\notin\{u_1,v_1\}$ and $x_3=u_2$.

\vspace{2mm}
Consider $G' = G^+ - x_2x_3$. 

\vspace{2mm}
If $G'$ is $2$-connected, then we consider a cycle $C^*$ traversing $x_1,z,x_2,v_2,x_3$ in this order (observe that $G'$ contains a cycle through $x_2$ and $x_3$, and $d_{G'}(z)=d_{G'}(x_2)=2$). Because of case (b) before we may assume that $x_1x_3\notin E(G')$. Therefore we denote $t_1=N_{C^*}(x_1)-\{z\}$, hence $t_1\notin \{v_2,x_3\}$. Now we apply Theorem B to obtain an $[x_1; t_1, v_2, x_3]$-$EPS$-graph $S^*=E^*\cup P^*$ of $G'$ with $C^*\subseteq E^*$. Now it is straightforward to see that $(S^*)^2$ has a hamiltonian cycle as required (containing an edge of $G$ incident to $x_i$ for both $i=1$ and $i=2$).

\vspace{2mm}
If $G'$ is not $2$-connected, we define $B_z$, $G_3$, and correspondingly $S_3$ as in   {\em Case (C)}(1)(b)(ii).

\vspace{2mm}
Let $C_z$ be a cycle in $B_z$ containing $z,x_1,x_2,v_2,t_1,c$, where $t_1=N_{C_z}(x_1)-\{z\}$.

\vspace{2mm}
If $c\notin \{x_1,v_2,t_1\}$, then by Theorem B let $S_z=E_z\cup P_z$ be an $[x_1; v_2, c, t_1]$-$EPS$-graph of $B_z$ with $C_z\subseteq E_z$. Now we continue as in {\em Case (C)}(1)(b)(ii)}, additionally using that $d_{P_z}\leq 1$.

\vspace{2mm}
If $c=v_2$, then again $v_2$ is a cutvertex of $G$, a contradiction.

\vspace{2mm}
Now suppose $c=x_1$. Hence $d_{G_3}(x_1)=1$ and $G_3$ is a non-trivial blockchain (note that $G_3=K_2$ is not possible because of $x_3\notin N(x_1)$ in this case). Now we continue as in the corresponding subcase of {\em Case (C)}(1)(b)(ii)} with both of $x_1$ and $x_2$ being pendant vertices of $G_2=B_z-z$.

\vspace{2mm}
Finally suppose $c = t_1$. Because of $C^+$ we conclude that $t_1=u_1$. We have $\{z, x_1, x_2, u_1, v_1, v_2\} \subset V(B_z)$. Thus there is a cycle $C_1 \subset B_z$  traversing $x_1u_1, x_1v_1$. Consequently, $z, x_2 \not \in C_1$. In fact, $ \widehat{C}  = C^+\triangle C_1 $  is a cycle containing $u_1, v_1,  x_1, z, x_2, x_3$ in this order; i.e., $x_1u_1$ is a chord of $\widehat{C}$. Hence $G'' = G^+ - x_1u_1$ is $2$-connected. Therefore $G''' = G'' - x_2x_3$ is a non-trivial block chain with one endblock $B_z''' \subset B_z$ since $G_3 \subset G''' $ ($z \in V(B_z''')$). Thus we can write

\vspace{2mm}  \hspace {10mm}
$G^+ - \{x_1u_1, x_2x_3\} = G''' = B_z''' \cup G_3'''$  \ \ \ with  \ \ \ $ B_z ''' \cap G_3''' = \{ c_0\} $,

\vspace{2mm}
\noindent where $c_0$ is a cutvertex of $G'''$ (possibly $c_0 = c$).

\vspace{1mm} 
Then $(G_3''')^2$ has a hamiltonian cycle $H_3$ containing $c_0w_0, x_3w_3 \in E(G_3''')$ by Corollary 1 {\it (i)} if $G_3'''$ is a non-trivial block chain, or by Theorem E if $G_3'''$ is 2-connected. Note that $G_3'''=K_2=u_1x_3$ is not possible because of $d_G(u_1)>2$ and $N(x_3)\subseteq V_2(G)$ in this case. 

\vspace{1mm} 
Let $C_z$ be a cycle in $B_z'''$ containing $z, x_1, x_2, c_0$.

\vspace{1mm} 
If $c_0 \not \in \{v_1, v_2\}$, we operate with a $[v_i; v_{3-i}, c_0]$-$EPS$-graph $S''' = E''' \cup P'''$ of $B_z'''$ with $C_z \subseteq E'''$, $i \in \{1, 2\}$, which exists by Theorem C. $(S''')^2$ contains a hamiltonian cycle $H'''$ containing $zx_1, zx_2, x_1v_1, c_0z_0, x_2v_2 \in E(B_z)$. $(H_3 - c_0w_0) \cup (H''' - c_0z_0) \cup \{w_0z_0\}$ is a required hamiltonian  cycle in $(G^+)^2$ with $x_1v_1, x_2v_2, x_3z_3 \in E(G)$.

\vspace{1mm} 
If $c_0 = v_i$ and $c_0 \neq v_{3-i}$, $i \in \{1, 2\}$, we operate with a $[v_i; v_{3-i}]$-$EPS$-graph $S''' = E''' \cup P'''$ of $B_z'''$ with $C_z \subseteq E'''$, which exists by Theorem D. $(S''')^2$ contains a hamiltonian cycle $H'''$ containing $zx_1,zx_2,x_1v_1,c_0z_0,x_2v_2 \in E(B_z)$. Again, $(H_3 - c_0w_0) \cup (H''' - c_0z_0) \cup \{w_0z_0\}$ is a required hamiltonian  cycle in $(G^+)^2$ with $x_1v_1, x_2v_2, x_3z_3 \in E(G)$.

\vspace{1mm} 
If $c_0 = v_1=v_2$, then $B_z''' = zx_1c_0x_2z$.  $(H_3 - c_0w_0) \cup (B_z''' - x_ic_0) \cup \{x_iw_0\}$ is a  required hamiltonian  cycle in $(G^+)^2$ with 
$x_{3-i}c_0, x_3z_3 \in E(G)$ for each $i \in \{1,2\}$.

\vspace{3mm}
(2)  Suppose $N(x_3) \not \subseteq V_2(G)$.

\vspace{2mm}
Then $d_G(x_3) = 2$. Set $N(x_3) = \{u_3, v_3\}$. Now we set $x_3^* \in N(x_3) - V_2(G)$. As before, let $C^+$ be of the  form $zx_1u_1 \cdots u_3x_3v_3 \cdots u_2x_2z$.

\vspace{2mm}
(a)  Suppose $x_3 \not \in N(x_1) \cup N(x_2)$.

\vspace{1mm} \hspace{5mm} (a1) $x_3^* \not \in \{u_1, u_2\}$.

\vspace{1mm} By Theorem B, let $S = E \cup P$ be an $[x_i; u_1, u_2, x_3^*]$-$EPS$-graph of $G^+$ with $C^+ \subseteq E$ for any $i \in \{1, 2\}$. Since $d_P(x_i) = 0 = d_P(z) $, $d_P(u_1) \leq 1$ and $d_P(u_2) \leq 1$, a required hamiltonian cycle $C_i$ in $(G^+)^2$ can be constructed in $S^2$ for each $i \in \{1, 2\}$ due to the restriction on $x_3^*$.

\vspace{1mm} \hspace{5mm} (a2) $x_3^*=u_3=u_1$ (the case $x_3^* = u_2 $ is symmetrical and therefore does not need separate consideration).

\vspace{1mm} (i) $v_3 \neq u_2$. In this case we operate with an $[x_i; u_1, v_3, u_2]$-$EPS$-graph $S_i = E_i \cup P_i$ with $C^+ \subseteq E_i$, $ i = 1, 2$ (see Theorem B). The restrictions on $x_3^*$ and $v_3$ guarantee that $(S_1)^2$ and $(S_2)^2$ yield hamiltonian cycles as claimed by the theorem.

\vspace{1mm} (ii) $v_3 = u_2$. That is, $C^+ = zx_1u_1x_3u_2x_2z$.

\vspace{1mm}
Assume first that one of $u_1, u_2$ is $2$-valent, i.e., $d_G(u_2) =2$ since $u_1=x_3^*\not\in V_2(G)$ . Then we operate with an $[x_i; u_1, x_{3-i}]$-$EPS$-graph $S_i = E_i \cup P_i$ with $C^+ \subseteq E_i$ for each $i \in \{1, 2\}$, which exists by Theorem~C. A required hamiltonian cycle in $(G^+)^2$ containing $x_iu_i, x_3u_1 \in E(G)$ can be constructed for each $i \in \{1, 2\}$.

\vspace{1mm}  Hence assume that $d_G(u_i) > 2$, $i=1, 2$. $G' = G-x_3$ is a trivial or non-trivial block chain.

\vspace{1mm} Suppose $G'$ is $2$-connected. Using induction,  $(G')^2$ has an $x_1x_2$-hamiltonian path $P_i' (x_1, x_2)$ containing $x_iw_i, u_iz_i \in E(G')$ for each $i \in \{1, 2\}$. Then $$ (P_i(x_1, x_2) - u_iz_i) \cup \{u_i x_3 z_i\}  $$ is a required hamiltonian path in $G^2$ for each $i \in \{1, 2\}$.

\vspace{1mm}
Finally assume that $G'$ is a non-trivial block chain. The endblock $B_i'$ in $G'$ containing $u_i, x_i$ is $2$-connected (since  $d_{G'}(u_i) \geq  2$), $i \in \{1, 2\}$; it also contains $v_i$ since $d_G(x_i) =2$. Let $P(x_i, u_i)$ denote an $x_iu_i$-path in $B_i'$ containing $v_i$ for any $i \in \{1, 2\}$.  Define the cycle
$\widetilde{C_i}$ by
$$ E(\widetilde{C_i}) = (E(C^+) -x_iu_i) \cup E(P(x_i, u_i)). $$   
Let $\widetilde{G_i} = G^+ - x_iu_i$, $i \in \{1, 2\}$; $\widetilde{G_i}$ is $2$-connected because $x_i u_i$ is a chord of $\widetilde{C_i}$.  $\widetilde{C_i}$  contains $z, x_1, x_2, x_3, u_1, u_2, v_i$. By Theorem C, there is a $[u_i; v_i, u_{3-i}]$-$EPS$-graph $\widetilde{S_i} = \widetilde{E_i} \cup \widetilde{P_i} \subset \widetilde{G_i}$ with $\widetilde{C_i} \subseteq \widetilde{E_i}$. A required hamiltonian cycle in $(\widetilde{S_i})^2$ containing $x_iv_i, x_3 u_i \in E(G)$ can be constructed, for each $i \in \{1, 2\}$.
\vspace{2mm}

(b) Suppose $x_3 \in N(x_1)$ but $x_3 \not \in N(x_2)$; that is, $x_3 =u_1$. By definition of $x_3^*$ we have $x_3^* = v_3$.

\vspace{1mm} Suppose $v_3 \neq u_2$, i.e., $N(x_3) \cap N(x_2) \cap C^+ = \emptyset$. To get a required hamiltonian cycle $C_i$ containing $x_iu_i, x_3x_3^*$, we operate with an $[x_i; u_2, x_3^*]$-$EPS$-graph $S_i= E_i \cup P_i$ of $G^+$ with $C^+ \subseteq E_i$, which exists by Theorem~C.

\vspace{1mm} Hence suppose $v_3 = u_2$, i.e., $N(x_3) \cap N(x_2) \cap C^+ \neq \emptyset$. Because of $d_G(x_3^*) > 2$ there exists $w_3 \in N(x_3^*) - \{x_3, x_2\}$. There is a $w_3v_2$-path $P(w_3, v_2)$ in $G$ not containing $x_3^*$ and therefore, $x_1,  x_2 \not \in P(w_3, v_2)$. Then $C^* = zx_1x_3u_2w_3P(w_3, v_2)v_2x_2z$ is a cycle in $G^+$ with $N(x_3) \cap N(x_2) \cap C^* = \emptyset$; thus we are back to the preceding case.

\vspace{2mm}

(c) Suppose $x_3 \not \in N(x_1)$ but $x_3  \in N(x_2)$.

\vspace{1mm}
This case is symmetrical to case (b) above.

\vspace{2mm}

(d) Suppose $x_3 \in N(x_1) \cap N(x_2)$. This case is not possible because of $N(x_i)\nsubseteq V_2(G)$, $i=1,2,3$.

\vspace{1mm}


\vspace{3mm}
(II)  Assume that $D(G) \neq \emptyset$.

\vspace{2mm}
We apply Theorem \ref{fletheorem1} to $G$ with respect to $\{x_1, x_2\}$ to conclude that $D(G)$ contains an edge $e$ such that $G-e$ has a $DT$-endblock $B_e$ such that $ A = \{x_1, x_2\} \not \subset V(B_e)$, and if $x_i\in V(B_e)$, then it is a cutvertex of $G-e$. Let $B_e'$ denote the other endblock of $G-e$. Also, let $c$ and $c'$ denote the cutvertices of $G-e$ belonging to $B_e$ and $B_e'$ respectively. If $c \neq c'$, set $G_0 = G-e -(B_e \cup B_e')$; it is a block chain containing $c, c'$ which are not cutvertices of $G_0$. Also, for the above $e$, denote $e=xx'$ where $x \in V(B_e)$ and $x' \in V(B_e')$.

\vspace{1mm}
Let $X = \{x_1, x_2, x_3\}$. Suppose $X \cap V(B_e) = \emptyset$. Then we replace the subgraph $B_e$ in $G$ with a path of length $3$ to obtain the $2$-connected edge-critical  graph $H$. By induction, $H$ has the strong ${\cal F}_3$ property. Moreover any strong ${\cal F}_3$ $x_1x_2$-hamiltonian path in $H^2$ can be converted into a strong ${\cal F}_3$ $x_1x_2$-hamiltonian path in $G^2$ by the method used in \cite{fle2:refer}. Hence we can assume that $X \cap V(B_e) \neq \emptyset$. We also note that it is tacitly assumed that the hamiltonian paths/cycles in the square of the smaller graphs contain as many edges of the given graphs as possible. The purpose of this assumption (already formulated in [7] and subsequent papers) is to facilitate the induction step and to keep the various cases arising, under control.

\vspace{1mm} 
With the  same argument as above, we see that $X \cap V(B_e') \neq \emptyset$ if $B_e'$ is a $DT$-block. If $B_e'$ is not a $DT$-block, then there is an edge $f \in E(B_e') \cap D(G)$ such that one of the endblocks $B_f$ of $G-f$ is a $DT$-block and $V(B_f) \subset V(B_e')$. This means that $X \cap V(B_f) = \emptyset$ if $X \cap V(B_e') = \emptyset$, and again the above argument  can similarly be applied. 

\vspace{1mm} 
In the ensuing discussion we keep in mind that there are at least two $DT$-endblocks $B^*$ and $B^{**}$ defined by the same element $e^*$ or by different elements $e^*, f^*\in D(G)$; and $B^*\cap B^{**}=\emptyset$, or $B^*\cap B^{**}=c^*$ where $c^*$ is a cutvertex of $G-e^*$. Therefore, a case not considered in $B^*$ implies a (sort of complementary) case in $B^{**}$ which is being taken care of when it occurs in $B^*$. 
	
\vspace{2mm}
Next we consider two special cases.

\vspace{3mm}
{\em Case (A):}  $X\cap (V(B_e)-c)=\emptyset$, or $x_3\in V(B_e)-\{c,x\}$ and $A\cap V(B_e)=\emptyset$. 

\vspace{1mm} In the first case it follows from the preceding argument that $c\in\{x_3,x_i\}$ for some $i\in\{1,2\}$. As before, we replace the subgraph $B_e$ in $G$ with a path of length $3$ to obtain the $2$-connected edge-critical  graph $H$. By induction, $H$ has the strong ${{\cal F}}_3$ property. Moreover, by a careful study of the method used in [7] one sees that any strong ${{\cal F}}_3$ $x_1x_2$-hamiltonian path $P_H$ in $H^2$ can be converted into a strong ${{\cal F}}_3$ $x_1x_2$-hamiltonian path in $G^2$. This applies, in particular, to the case where $c\in A$ and $P_H$ contains an edge of $H$ incident to $c$ (here, some of the 13 cases listed in [7] need not be considered). Hence we are left with the case where $X\cap V(B_e)=x_3$ and $x_3\not= c,x$.

We proceed as before, replacing $B_e$ with a path $P_3$ of length 3; again, the resulting graph is denoted by $H$. By induction on $|D(G)|,H^2$ has a strong ${{\cal F}}_3$ $x_1x_2$-hamiltonian path containing $e_u\in E(G)$ incident to $u$ for some $u\in V(G)-V(B_e)$. In fact, a careful study of the procedure employed before shows that $P_H$ can be converted into a strong ${\cal F}_3$ $x_1x_2$-hamiltonian path $P_{1,2}$ of $G^2$ containing an edge of $G$ incident to $x_3$. Namely, depending on the various cases of the traversal of $V(P_3)$ by $P_H$, 
\begin{itemize}
	\item 
one either applies Lemma 3 to use a hamiltonian cycle $C_e$ of $(B_e)^2$ such that $C_e$ traverses in $c$ edges of $B_e$, and likewise, $C_e$ traverses at least one edge in $x$ and at least one edge in $x_3$, belonging to $B_e$ (observe that $|V(B_e)|\geq 4$ since $G$ is edge-critical and thus does not have a triangle); 
\item
or one applies induction to use a hamiltonian path $P_e$ of $(B_e)^2$ joining $x$ and $c$ and containing at least one edge of $B_e$ incident to $x_3$ and an edge of $B_e$ incident to any given $t\in\{x,c\}$.  
\end{itemize}
\vspace{1mm}

{\em Case (B):} $X\cap V(B_e)=x_3=x$.

It follows that $A\cap V(B'_e)\not =\emptyset$; without loss of generality $x_1\in V(B'_e)$. Assume the notation chosen in such a way that $x_1\not= c'$ if $A\subset V(B'_e)$. Moreover, if $A\cap V(B'_e)=x_1$ it follows from the preceding considerations that $x_1$ belongs to the DT-endblock $B_f\subseteq B'_e$ for some $f\in D(G)$ if $E(B'_e)\cap D(G)\not=\emptyset$, or else $B_{e'}$ is a DT-endblock; and $x_1\not=c_f$ by Theorem G, where $c_f$ is the cutvertex of $G-f$ in $B_f$. Also, $c_f=c'$ if $c'\in V(B_f)$. Hence $x_1\not= c'$ can be assumed in any case.

Denote the blocks of $G-e$ by $B_0,\dots,B_k$ according to their order in $bc(G-e)$ such that $B_0=B'_e$, $B_k=B_e$, and let $j$ be the smallest index such that $x_2\in V(B_j)$; possibly $j=0$. Set 
$$G_{0,j}:=\bigcup^j_{i=0}B_i\mbox{ \rm and }G_{j+1,k}:=\bigcup^k_{i=j+1}B_i.$$ 

(i) Suppose $j>0$. By applying induction to the individual $2$-connected blocks of $G_{0,j}$ it follows that $(G_{0,j})^2$ has a hamiltonian $x_1x_2$-path $P_{1,2}$ containing $x_iy_i\in E(G_{0,j}), i=1,2$, for some $y_1,y_2\in V(G_{0,j})$, as well as $x'y'\in E(B'_e), c^*y^*\in E(B_j)$ where $c^*=B_j\cap B_{j+1}$, and $x'y'=x_1y_1$ if $x'=x_1$, $c^*y^*=x_2y_2$ if $c^*=x_2$. Likewise by Corollary 1{\it (i)}, $(G_{j+1,k})^2$ has a hamiltonian cycle $\tilde C$ containing $x_3y_3$, $x_3z_3\in E(B_e)$ and $c^*z^*\in E(B_{j+1})$. 

In any case, 
$$(P_{1,2}\cup\tilde C-\{c^*y^*, c^*z^*\})\cup\{ y^*z^*\} \hskip 2.0cm (1)$$
defines a hamiltonian $x_1x_2$-path of $G^2$ containing $x_1y_1$ and $x_3y_3$; it also contains $x_2y_2$ if $x_2\not=c^*$. On the other hand, if $c^*=x_2$ we construct a hamiltonian $x_1x_2$-path of $G^2$ containing $x_2y_2$ and $x_3y_3$ as follows: if $j+1<k$ and $\kappa(B_{j+1})\geq 2$, then $\tilde C$ can be assumed to contain $c^*z^*$, $c^*x^*\in E(B_{j+1})$, and we set $$\tilde{\tilde C} =(\tilde C-\{c^*z^*,c^*x^*\})\cup\{x^*z^*\} $$  which defines a hamiltonian cycle of $(G_{j+1,k}-c^*)^2$ containing $x_3y_3$, $x_3z_3$. The same type of hamiltonian cycle is obtained if $B_{j+1}$ is a bridge of $G-e$.  Thus, in both cases
$$(P_{1,2}\cup\tilde{\tilde C}-\{x'y', x_3z_3\})\cup\{x'z_3,y'x_3\}\hskip 2.0cm (2)$$
defines a hamiltonian $x_1x_2$-path of $G^2$ containing $x_2y_2$ and $x_3y_3$. 

Thus we are left with the case $j+1=k$ implying $x_2\not=c$ and thus $\kappa(B_j)\geq 2$. We now proceed as in (1) above. 

(ii) $j=0$. That is, $G_{0,0}=B'_e$; $\{x_1,x_2\}\subseteq V(B'_e)$ follows. 

If $\{x_1,x_2\} \not=\{x',c'\}$ we obtain by induction a hamiltonian $x_1x_2$-path $(P^{(i)}_{1,2})$  of $(B'_e)^2$ containing $x_iy_i\in E(B'_e)$ for any $i\in\{1,2\}$, but also an edge $e_{t'}$ incident to $t'\in\{x',c'\} - \{x_1,x_2\}$. If $t'=c'$ we proceed as in (1), whereas we proceed as in (2) if $t'=x'$. 

Thus we assume $\{x_1,x_2\}=\{x',c'\}$; by the initial choice of notation, $x_1=x'$ and $x_2=c'$ follows. Also, $c'\not=c$ by the hypothesis of this {\it Case (B)}.

Since  $c'\not=c$, $G':=G-e-B'_e$ is a non-trivial block chain. By Corollary 1{\it (ii)}, $(G')^2$ has a hamiltonian $x_2x_3$-path $P_{2,3}$ containing edges $x_2y_2, x_3y_3\in E(G')$. By Theorem E, $(B'_e)^2$ has an $[x_i;x_{3-i}]$-hamiltonian cycle $C_i$ for every $i\in\{1,2\}$. Denote the corresponding edges of $C_i\cap E(B'_e)$ by $x_1y_1^{(1)}, x_1z_1^{(1)}$ and $x_2z_2^{(1)}$ if $i=1$, and by $x_1y_1^{(2)}$ and $x_2u_2^{(2)}, x_2z_2^{(2)}$ if $i=2$.

Assume further the notation chosen in such a way that the $x_1x_2$-path in $C_i$ containing $x_1y_1^{(i)}$ also contains $x_2z_2^{(i)}$; denote it by $P_1^{(i)}$ and set $P_2^{(i)}=C_i-P_1^{(i)}$.
$$(P_1^{(i)}-x_2)\cup\{z^{(1)}_2y_2\}\cup(P_{2,3}-x_2)\cup\{x_3z_1^{(1)}\}\cup(P_2^{(1)}-x_1)$$
defines a hamiltonian $x_1x_2$-path of $G^2$ starting with $x_1y_1^{(1)}\in E(G)$. Likewise
$$(P_2^{(2)}-x_2)\cup\{u_2^{(2)}z_2^{(2)}\}\cup(P_1^{(2)}-\{x_1,x_2\})\cup\{y_1^{(2)}x_3\}\cup P_{2,3}$$
defines a hamiltonian $x_1x_2$-path of $G^2$ ending with $x_2y_2\in E(G)$. 

This finishes {\it Case (B)}.

\vspace{3mm}

For the remaining cases of the proof of Theorem 3 we consider $D_1(G)$ comprising those edges of $D(G)$ such that for every $g\in D_1(G)$ one of the endblocks of $G-g, $ $B_g$ say, is a $DT$-graph. Let $c_g$ denote the cutvertex of $G-g$ in $B_g$. Having solved the {\it Cases (A)} and {\it (B)} we conclude that $X\cap(V(B_g)-c_g)\not=\emptyset$ in any case. Note that $G$ has at least two $DT$-endblocks as just described: one is the aforementioned $B_e$, another one is either $B'_e$, or $B'_e$ contains $f\in D_1(G)$ such that the corresponding $DT$-endblock $B_f$ is a proper subgraph of $B'_e$. $B_e\cap B_f=\emptyset$ or $B_e\cap B_f=c$ where $c=c_f=c_e$.
	
We proceed analogous to {\it Case (B)} denoting the blocks of $G-e$ by $B_0,\dots, B_m$ with $B_0=B_e$, $B_m=B'_e$.
	
Suppose first that $|X\cap V(B_e)|=1$. In view of {\it Cases (A)} and {\it (B)} we have 
	$$X\cap V(B_e)=A\cap V(B_e)=A\cap(V(B_e)-c)=x_i  \ ;$$ 
without loss of generality $i=1$. By the same token $x_2\in V(B_f)\subseteq V(B_e')$ and $x_2\neq c'$.

Let $P(c,x_1) $ be a $c,x_1$-hamiltonian path of $B^2_e$ with $x_1z_1\in(P(c,x_1))\cap E(B_e)$, which exists by Theorem F {\it (ii)}.	
	
If $x_3\in V(B_e')-\{c'\}$, we operate with an $x_2c'$-hamiltonian path $P(x_2,c')$ of $(B'_e)^2$ with $x_2z_2, x_3z_3 \in E(P(x_2,c'))\cap E(B'_e)$ using induction, and trivially with a $c'c$-hamiltonian path $P_0$ of $G^2_0$. Note that $P_0=\emptyset$ if $c=c'$ in this case.

If $x_3\notin V(B_e')-\{c'\}$, we operate with an $x_2c'$-hamiltonian path $P(x_2,c')$ of $(B'_e)^2$ with $x_2 \in E(P(x_2,c'))\cap E(B'_e)$ which exists by Theorem F {\it (ii)}, and with a $c'c$-hamiltonian path $P_0$ of $G^2_0$ containing $x_3z_3\in E(G_0)$ applying Theorem F to each 2-connected block of $G_0$. Note that $c=c'=x_3$ is not possible by the assumption $X\cap V(B_e)=x_1$ and it covers the case $x_3=c'\neq c$.

Then 
	$$P(x_2,c'), P_0,P(c,x_1) $$
is a hamiltonian $x_1x_2$-path of $G^2$ containing $x_3z_3,x_iz_i\in E(G)$, for $i=1,2$. This settles the case $|X\cap V(B_e)|=1$.
	
Now suppose $|X\cap V(B_e)|=2$. Because of the case just settled we must also have $|X\cap V(B_f)|=2$ implying $c_f=c'=c\in X$. Again, suppose without loss of generality that $x_1\in B_e$.
	
If $c=x_3$, let $P_0$ be an $x_1x_3$-hamiltonian path of $B^2_e$ such that $x_1z_1\in E(P_0)\cap E(B_e)$ and $P_{1}$ be an $x_3x_2$-hamiltonian path of $B^2_{e'}$ such that $x_3z_3\in E(P_1)\cap E(B_{e'})$, which exist by Theorem F {\it (ii)}. Hence 
	$$P_0P_1$$
is an $x_1x_2$-hamiltonian path of $G^2$ with $x_1z_1, x_3z_3\in E(G)$. We proceed analogously to obtain an $x_2x_1$-hamiltonian path of $G^2$ with $x_2z_2, x_3z_3 \in E(G)$ as required by the theorem.
	
Now suppose without loss of generality that $c=x_1$ and $x_3\in V(B_e)-c$. Hence we have $x_1,x_2\in V(B_f)$, i.e., $x_1,x_2\in V(B'_e)$. We apply Theorem E to $B_e$ and either Theorem~F{\it (ii)} or Theorem 2 to $B'_e$.
	
By Theorem E, $B^2_e$ contains a hamiltonian cycle $C_e$ with $x_1y_1,x_1z_1\in E(C_e)\cap E(B_e)$ and $x_3u_3\in E(C_e)\cap E(B_e)$. 
	
As for $(B'_e)^2$, it has an $x_1x_2$-hamiltonian path $P_{1,2}$ with $x_1u_1\in E(P_{1,2})\cap E(B'_e)$, by Theorem~F{\it (ii)}. Thus,
$$E(P_{1,2})\cup E(C_e)\cup \{u_1y_1\}-\{x_1u_1, x_1y_1\}$$
defines on $x_1x_2$-hamiltonian path of $G^2$, containing $x_1z_1\in E(G)$, but also $x_3u_3\in E(G)$.
	
Likewise, Theorem 2 implies that $(B'_e)^2$ has either an $x_1x_2$-hamiltonian path $P_{1,2}$ with $x_1u_1,x_2z_2\in E(P_{1,2})\cap E(B'_e)$, or it has an $x_1x_2$-hamiltonian path $P_{1,2}$ with $x_2z_2\in E(P_{1,2})\cap E(B'_e)$ and $u_1v_1\in E(P_{1,2})$ for some $u_1, v_1\in N(x_1)$. In the first case, we define an $x_1x_2$-hamiltonian path of $G^2$ as above; it contains $x_2z_2, x_3u_3\in E(G)$. In the second case we proceed similarly: here,
	$$E(P_{1,2})\cup E(C_e)\cup\{u_1y_1,v_1z_1\}-\{u_1v_1,x_1y_1,x_1z_1\}$$
	defines an $x_1x_2$-hamiltonian path containing $x_2z_2,x_3u_3\in E(G)$. Thus $G$ has the strong ${\cal F}_3$ property.
	
	The case $X\subset V(B_e)$ needs no separate consideration since it implies $|X\cap V(B_f)|\leq 1$, in which case we may consider $B_f$ instead of $B_e$. This finishes the proof of  Theorem 3.
\vspace{1mm}

\section{Arbitrary $2$-connected graphs}

We now proceed to prove the main result of this paper.

\vspace{2mm}
\begin{result} \label{f4}
Let $G$ be a $2$-connected graph.  Then $G$ has the ${\cal F} _4$ property.
\end{result}

\vspace{2mm}
\noindent
{\bf Proof:} \ We may assume that $G$ is an edge-critical block since otherwise we can delete edges of $G$ until we reach an edge-critical block.

\vspace{1mm} If $G$ is a $DT$-block, then the result is true by Theorem \ref{dt}. So assume that $G$ is not a $DT$-block. The rest of the proof is by induction on $|D(G)|$, or on $|V(G)|$. That is, if $H$ is an edge-critical block with $|D(H)| < |D(G)|$ or $|V(H)| < |V(G)|$, then $H$ has the ${\cal F} _4$ property.

\vspace{2mm}
   By \cite[Theorem 1]{fle2:refer},  $D(G)$ contains an edge $e$ such that $G-e$ is a block chain  with  at least one of its endblocks, say $B_e$, being a $DT$-block.  Let $B_e'$ be the other endblock of $G-e$.

\vspace{2mm} 
Throughout, we let $e=xx'$ where $x \in V(B_e)$ and $x' \in V(B_e')$.

\vspace{1mm}
We claim that $D(G)$ contains an edge $e^*$ such that $G - e^*$ has an endblock $B_{e^*}$ which is a $DT$-block satisfying $|V(B_e)\cap V(B_{e^*})| \leq 1$. To see this, we note that if $B_e'$ is also a $DT$-block, then $e^* = e$ and $B_e' = B_{e^*}$, and the inequality holds trivially. If $B_e'$ is not a $DT$-block, then it is edge-critical and again  \cite[Theorem 1]{fle2:refer} applies and  $e^* $ is in $D(G) \cap B_e'$, and $B_{e^*}$ is a subgraph of $B_e'$. Since $|V(B_e)\cap V(B'_e)|\leq 1$ the claimed inequality holds.

\vspace{1mm} Let $X = \{x_1, x_2, x_3, x_4\}$ and let $k = \ {\rm min} \ \{ |V(B_e) \cap X|, \  |V(B_{e^*}) \cap X|\}$. Then clearly $k \leq 2$. 
Without loss of generality, we assume that $|V(B_e) \cap X|=k$.

\vspace{2mm}
We first dispose of the case $k=0$ proceeding as in the proof of Theorem 3: we replace $B_e$ by a path of length $3$. The resulting graph $H$ is an edge-critical block and $|D(H)| < |D(G)|$. By induction $H$ has the ${\cal F} _4$ property. Any ${\cal F} _4$ hamiltonian path in $H^2$ can then be converted into an ${\cal F} _4$ hamiltonian path in $G^2$ by the same method  used  in \cite{fle2:refer}.

\vspace{2mm}
Hence we assume that $k \in \{1, 2\}$.

\vspace{2mm}
Let $c$ and $c'$ be the cutvertices of $G-e$ belonging to $B_e$ and $B_e'$ respectively.   Note that if  $c=c'$, then $G-e$ is a block chain with only $2$ blocks $B_e$ and $B'_e$.

\vspace{2mm}   If $k=2$, then we may assume without loss of generality  that either $V(B_e) \cap X = \{x_3,  x_4\}$ or $V(B_e) \cap X = \{x_2,  x_4\}$, or $V(B_e)\cap X=\{x_1,x_2\}$ (namely, if $c=c'=x_i\in\{x_1,x_2\}$, $\{x_3,x_4\}\subseteq V(B_{e^*})$, and $B_e'=B_{e^*}$).

\vspace{2mm} If $k=1$, then we may  assume without loss of generality  that either $V(B_e) \cap X = \{x_2\}$ or $V(B_e) \cap X = \{x_4\}$. In any case, we note that $c  \not \in \{x_2, x_4\}$: otherwise, we replace $B_e$ by a path of length $3$ to obtain $H$ which has an ${\cal F} _4$ $x_1x_2$-hamiltonian path in $H^2$. Again, as before, we apply the method used in \cite{fle2:refer} to see that  any corresponding ${\cal F} _4$ hamiltonian path in $H^2$ can be converted into an ${\cal F} _4$ hamiltonian path in $G^2$.

\vspace{2mm}
{\em Case (A):} $c=c'$

\vspace{2mm} (1) Suppose $k=2$.

\vspace{2mm} (a) Suppose $x_3, x_4 \in V(B_e)$.  Then $x_1, x_2 \in V(B_e')$ and there is an  $x_1x_2$-hamiltonian path  $P'(x_1,x_2)$ in $(B_e')^2$ containing  an edge $cw'$ of $B_e'$, by  Theorem \ref{fletheorem4} (note that $c  \not \in \{x_1, x_2\}$ by assumption). Let $w \in N(c) \cap V(B_e)$. Then, by Theorem \ref{dt}, there is an ${\cal F} _4$ $cw$-hamiltonian path  $P(c,w)$ in $(B_e)^2$ containing $x_3z_3, x_4z_4$ which are edges of $B_e$ if $\{c, w\} \cap \{x_3, x_4\} = \emptyset$.

\vspace{1mm} Suppose  $|\{c, w\} \cap \{x_3, x_4\}| =1$.  Without loss of generality,  assume that $x_3 \in \{c, w\}$.  By Theorem \ref{strongf3}, $B_e$ has the strong ${\cal F}_3$ property. Consequently, $B_e^2$ has a $cw$-hamiltonian path $P(c,w)$ containing $x_4z_4$, $cz\in E(B_e)$, $cz\not=cw$ if $x_3=c$; or it contains $x_4z_4$, $wv\in E(B_e)$, $wv\not=cw$, if $x_3=w$.

 \vspace{1mm} Suppose  $\{c, w\} = \{x_3, x_4\}$. If $d_{B_e}(c)>2$, then consider $u\in N(c)\cap V(B_e)-w$ such that $u\not\in\{x_3,x_4\}$ and argue with $u$ in place of $w$ as in the preceding case. Thus we may assume that $d_{B_e}(c)=2$. By Theorem \ref{fletheorem3},  $(B_e)^2$ has a $[c;w]$-hamiltonian cycle $C_w$ containing $cw, vw, cz$ which are three different  edges of $B_e$.  Let $C_w - cw = P(c,w)$.

 \vspace{1mm}
By deleting in all cases the edge $cw'$ from $P'(x_1, x_2)$ and adding $ww' \in E(G^2)$, we have a required ${\cal F} _4$ $x_1x_2$-hamiltonian path in $G^2$.

\vspace{2mm} (b) Suppose $x_2, x_4 \in V(B_e)$.  Then $x_1, x_3 \in V(B_e')$. 
Note that $c \not \in \{x_1, x_3\}$ since $k =2$.

\vspace{2mm} 
If $x_2 \neq c$, then by Theorem \ref{fletheorem4}, there is an  $x_2c$-hamiltonian path  $P(x_2,c)$ in $(B_e)^2$ containing  an edge $x_4z_4$ of $B_e$ (independent  of $x_4=c$ or $x_4 \neq c$)  and there is an  $x_1c$-hamiltonian path  $P'(x_1,c)$ in $(B_e')^2$ containing   $x_3z_3 \in E(B'_e)$. $P'(x_1,c)$ and $P(x_2,c)$ form a required ${\cal F} _4$ $x_1x_2$-hamiltonian path in $G^2$.

\vspace{2mm} If $x_2=c$, we apply Theorem \ref{fletheorem3} to $B_e$  to obtain a hamiltonian cycle $C$ in $(B_e)^2$ containing $x_2v_1, x_2v_2, x_4z_4$ which are edges of $B_e$.  By Theorem \ref{strongf3}, $(B_e')^2$ has a strong ${\cal F}_3$ $x_1x_2$-hamiltonian path $P'(x_1, x_2)$ containing $x_2w', x_3z_3$ which are edges of $B_e'$. A required ${\cal F} _4$ $x_1x_2$-hamiltonian path in $G^2$ is given by $(C-x_2v_1) \cup (P'(x_1, x_2)- x_2w') \cup \{v_1w'\}$.

\vspace{2mm} (c) Suppose $V(B_e)\cap X=\{x_1,x_2\}$ with $x_i=c$ for some $i\in\{1,2\}$; without loss of generality $x_2=c=c'$. $x_3,x_4\in V(B'_e)$ follows (note that the case $c\not\in\{x_1,x_2\}$ can be treated symmetrically to case (a)).  By Theorem \ref{fletheorem4}, $(B_e)^2$ has an $x_1x_2$-hamiltonian path $P(x_1,x_2)$ containing an edge $x_2v\in E(B_e)$.

\vspace{2mm} Suppose $w'\in (N(c')-\{x_3,x_4\})\cap V(B'_e)$ exists. By induction, $(B'_e)^2$ has a $w'x_2$-hamiltonian path $P'$ containing edges $x_3z_3,x_4z_4\in E(B'_e)$. Clearly,
$$(P(x_1,x_2)\cup P'-x_2v)\cup\{w'v\}$$
defines an $x_1x_2$-hamiltonian path of $G^2$ as required.

Finally suppose $N(c')\cap V(B'_e)=\{x_3,x_4\}$. That is, $c'=c=x_2$ is $2$-valent in $B'_e$. In this case we apply Theorem \ref{fletheorem3} to obtain a $[c';x_3]$-hamiltonian cycle  $C'$ in $(B'_e)^2$ containing three different edges $c'x_3,c'x_4,x_3z_3\in E(B'_e)$.

It follows that 
$$(P(x_1,x_2)-x_2v)\cup (C'-x_3c')\cup \{vx_3\}$$
defines an ${\cal F}_4$ $x_1y_2$-hamiltonian path of $G^2$ as required. This settles case (1).

\vspace{2mm} (2) Suppose $k=1$.

\vspace{2mm} (a) Suppose $x_2 \in V(B_e)$.  Then $x_1, x_3, x_4 \in V(B_e')-c$. Hence by induction there is an ${\cal F} _4$ $x_1c$-hamiltonian path  $P'(x_1,c)$ in $(B_e')^2$ containing $x_3z_3, x_4z_4\in E(B_e')$. In $(B_e)^2$, there is an $x_2c$-hamiltonian path  $P(x_2,c)$  which together with $P'(x_1,c)$ form a required ${\cal F} _4$ $x_1x_2$-hamiltonian path in $G^2$.

\vspace{2mm}
(b) Suppose $x_4 \in V(B_e)$.  Then $x_1, x_2, x_3 \in V(B_e')-c$ and by induction there is an ${\cal F} _4$ $x_1x_2$-hamiltonian path  $P'(x_1,x_2)$ in $(B'_e)^2$ containing $x_3z_3, cw'\in E (B'_e)$. In $(B_e)^2$, there is a hamiltonian cycle $C_c$ containing three different edges $cw, cz, x_4z_4\in E(B_e)$ by Theorem E. Delete $cw'$ from $P'(x_1, x_2)$ and $cw$ from $C_c$ and  join $w'$ to $w$ to obtain a required ${\cal F} _4$ $x_1x_2$-hamiltonian path in $G^2$. This settles case (2) and thus finishes the proof of {\it Case (A)}.

\vspace{2mm}
{\em Case (B):} $c \neq c'$

\vspace{2mm}
Let $G_0 = G -(e\cup B_e \cup  B_e')$.

\vspace{2mm} (1) Suppose $k=2$.

\vspace{2mm} In this case  $(V(G_0)-\{c,c'\}) \cap X = \emptyset$.  By Corollary \ref{flecor}, $(G_0)^2$ has a hamiltonian cycle $C_0$ containing $c'w_0', cw_0$ which are edges of $G_0$, provided $G_0$ is a non-trivial block chain. If, however, $G_0\not=K_2$ is a block, then such hamiltonian cycle $C_0$ exists  by Theorem \ref{fletheorem3}. Moreover, we only have to deal with the cases (1.1), (1.2) below; otherwise we could consider $B_e* \subseteq\ B_e'$.

\vspace{2mm} \hspace{5mm} (1.1) Suppose $x_3, x_4 \in V(B_e)$.

\vspace{1mm} Then $x_1, x_2 \in V(B_e')$.  If  $c' \not \in \{x_1, x_2\}$, then by Theorem \ref{fletheorem4}{\it (i)}, $(B_e')^2$ has an ${\cal F} _3$ $x_1x_2$-hamiltonian path $P'(x_1, x_2)$  containing an edge $c'w'$ of $B_e'$. If $c' \in \{x_1, x_2\}$, say $c'=x_1$, then we let $P'(x_1, x_2)$ denote an $x_1x_2$-hamiltonian path in $(B_e')^2$ containing an edge $x_1w' =c'w'$ of $B'_e$ (see Theorem \ref{fletheorem4}(ii)).

\vspace{1mm} (a) Suppose $c = x_i$ for some $i \in \{3, 4\}$. Let $C_c$ denote an $[x_i ; x_{7-i}]$-hamiltonian cycle in $(B_e)^2$ containing $x_iz_i, x_iw_i, x_{7-i}z_{7-i}$ which are edges of $B_e$. In this case,
 \[(P'(x_1, x_2)-c'w') \cup (C_0 - \{c'w_0', cw_0\}) \cup (C_c - x_iw_i) \cup \{w'w_0', w_0w_i\}\]
defines a required ${\cal F} _4$ $x_1x_2$-hamiltonian path in $G^2$ provided $G_0\not=K_2$. If, however, $G_0=K_2$, then 
$$(P'(x_1,x_2)-c'w')\cup (C_c-x_iw_i)\cup \{c'w_i,cw'\}$$
yields the required result.

\vspace{1mm} (b) Suppose $c \neq x_i$ for any $i \in \{3, 4\}$.

\vspace{1mm} (i) Suppose $N(c) \cap \{x_3, x_4\}= \emptyset$. Let $w \in N(c) \cap V(B_e)$.  By induction, there is an ${\cal F} _4$ $cw$-hamiltonian path $P_e(c, w)$ in $(B_e)^2$ containing $x_3z_3, x_4z_4$ which are edges of $B_e$. In this case,
\[(P'(x_1, x_2)-c'w') \cup (C_0 - \{c'w'_0, cw_0\}) \cup P_e(c, w) \cup \{w'w_0', ww_0\}\]
yields a required ${\cal F} _4$ $x_1x_2$-hamiltonian path in $G^2$; and if $G_0 = K_2$, then we obtain the required result analogously as in case (a).

\vspace{1mm} (ii) Hence we assume that  $N(c) \cap \{x_3, x_4\} \neq \emptyset$.

\vspace{1mm} If there exists $w \in N(c) \cap V(B_e) $ such that  $w \not \in \{x_3, x_4\}$, then the argument used in (i) applies and we have a required  ${\cal F} _4$ $x_1x_2$-hamiltonian path in $G^2$ as before.

\vspace{1mm}
So assume that $N(c) \cap V(B_e) =  \{x_3, x_4\} $. Let $C_c$ denote an $[c; x_3]$-hamiltonian cycle in $(B_e)^2$ containing $x_3c, x_3z_3, x_4z_4$ which are edges of $B_e$. Then
\[(P'(x_1, x_2)-c'w') \cup (C_0 - \{c'w_0', cw_0\}) \cup (C_c - x_3c) \cup \{w'w_0', w_0x_3\}\]
yields a required ${\cal F} _4$ $x_1x_2$-hamiltonian path in $G^2$ if $G_0\not=K_2$; and the case $G_0=K_2$ is treated analogously as before. 

\vspace{2mm} \hspace{5mm} (1.2) Suppose $x_2, x_4 \in V(B_e)$.

\vspace{1mm} Then $x_1, x_3 \in V(B_e')$.  If  $x_1 \neq c' $, then by Theorem \ref{fletheorem4}, $(B_e')^2$ has an ${\cal F} _3$ $x_1c'$-hamiltonian path $P'(x_1, c')$  containing an edge $x_3z_3$ of $B_e'$ (even if $x_3 = c'$). If $c' = x_1$, then by Theorem \ref{fletheorem3}, $(B_e')^2$ has an $[x_1; x_3]$-hamiltonian cycle $C'$ containing three edges $x_1w_1, x_1z_1, x_3z_3\in E(B'_e)$. In this case, let $P'(x_1, w_1) = C' - x_1w_1$.

\vspace{1mm}
Consider $B_e$. If  $x_2 \neq c $, then by Theorem \ref{fletheorem4}, $(B_e)^2$ has an ${\cal F} _3$ $x_2c$-hamiltonian path $P(x_2, c)$  containing an edge $x_4z_4$ of $B_e$ (even if $x_4 = c$). If $c = x_2$, then by Theorem \ref{fletheorem3}, $(B_e)^2$ has an $[x_2; x_4]$-hamiltonian cycle $C$ containing $x_2w_2, x_2z_2, x_4z_4$ which are edges of $B_e$. In this case, let $P(x_2, w_2) = C - x_2w_2$.

\vspace{1mm}
(a) Suppose $G_0\not=K_2$.

\vspace{1mm} By Corollary \ref{flecor}{\it (ii)}, Theorem F respectively, $(G_0)^2$ has a $cc'$-hamiltonian path $P_0(c, c')$ containing an edge $cw_0$  of $G_0$ incident to $c$, or an edge  $c'w_0'$ of $G_0$ incident to $c'$.  In the case that $G_0$ has $2$ or more blocks, then $P_0(c, c')$ can be chosen to contain both  $cw_0$ and $c'w_0'$.

\vspace{1mm}
(i)  Suppose $c\neq x_2$ and $c' \neq x_1$. Then  $P'(x_1, c') \cup P_0(c', c) \cup P(x_2, c) $  yields a required ${\cal F} _4$ $x_1x_2$-hamiltonian path in $G^2$.

\vspace{1mm}
(ii) Suppose $c \neq x_2$ and $c' = x_1$. Then $P'(x_1, w_1) \cup (P_0(c', c) - c'w_0') \cup P(x_2, c) \cup \{w_0'w_1\} $   yields a required ${\cal F} _4$ $x_1x_2$-hamiltonian path in $G^2$.

\vspace{1mm}
(iii) Suppose $c = x_2$ and $c' \neq  x_1$. Then $P'(x_1, c') \cup (P_0(c', c) - cw_0) \cup P(x_2, w_2) \cup \{w_0w_2\} $     yields a required ${\cal F} _4$ $x_1x_2$-hamiltonian path in $G^2$.

\vspace{1mm}
(iv) Suppose $c = x_2$ and $c' = x_1$.

\vspace{1mm} First assume that $G_0$ has $2$ or more blocks.
Then $$P'(x_1, w_1) \cup (P_0(c', c) - \{c'w_0', cw_0\}) \cup P(x_2, w_2) \cup \{w_0'w_1, w_0w_2\} $$ yields a required ${\cal F} _4$ $x_1x_2$-hamiltonian path in $G^2$.

\vspace{1mm} Next assume that $G_0$ is  $2$-connected.  By Theorem \ref{lemmaf3}, $(G_0)^2$ has a $cc'$-hamiltonian path $P_0(c,c')$  containing an edge $cw_0$ of $G_0$ and $P_0(c, c')$  either contains an edge $c'w_0'$ of $G_0$    or else contains an edge $uv$ for some vertices $u, v \in N(c') \cap V(G_0)$.   In the former case we proceed as in the case where $G_0$ has $2$ or more blocks,
to obtain a required $x_1x_2$-hamiltonian path in $G^2$. In the latter case,
$$P(x_2,w_2) \cup (P_0(c,c') -\{cw_0, uv\}) \cup (P'(x_1,w_1)-\{z_1x_1\}) \cup \{w_2w_0, w_1v, z_1u\}$$ 
yields a required $x_1x_2$-hamiltonian path in $G^2$.

\vspace{2mm}
(b) Suppose $G_0=K_2$.

\vspace{1mm} If $c \neq x_2$ or $c' \neq x_1$, then  the methods used in the above cases (a) (i), (ii), (iii) can be used to construct  a required ${\cal F} _4$ $x_1x_2$-hamiltonian path in $G^2$. Hence we assume that $c=x_2, c' = x_1$. Then by Theorem \ref{fletheorem4}, $(B_e)^2$ (respectively $(B_e')^2$) has an ${\cal F} _3$ $x_2x$-hamiltonian path   $P(x_2, x)$ containing $x_4z_4$ (respectively $x_1x'$-hamiltonian path $P'(x_1, x')$ containing $x_3z_3$) where $x_3z_3, x_4z_4 \in E(G)$ (even if $x_4=x$ and $x_3=x'$). Then  $P(x_2,x) \cup \{xx'\} \cup P'(x_1, x')$ is a required ${\cal F} _4$ $x_1x_2$-hamiltonian path in $G^2$.

\vspace{2mm} (2) Suppose $k=1$.

\vspace{2mm} Recall that, in this case,  either $V(B_e) \cap X = \{x_2\}$ or else  $V(B_e) \cap X = \{x_4\}$  and that $c \not \in \{x_2, x_4\}$.

\vspace{2mm}   (2.1) Suppose $x_2 \in V(B_e)$ and $x_2 = x$.

\vspace{2mm} Write the block chain $G-e$ as $B_1 \cup B_2 \cup \cdots \cup B_k$, $k >2$ with $B_i \cap B_{i+1} = c_i$ for $i=1, 2, \ldots, k-1$ where $B_1 = B_e'$ and $B_k = B_e$ so that $c_1 = c'$, $c_{k-1} = c$ and $G_0 = B_2 \cup \cdots \cup B_{k-1}$.

\vspace{2mm}
(a) Suppose $(V(B_e') - c_1) \cap X = \emptyset$ and $c_1=x_1$.

\vspace{1mm}
Then either (i) $B_e'$ is a $DT$-graph  or else (ii) $B_e'$ contains an edge $f \in D(G)$ such that one of its endblocks $B_f$ of $G-f$ is a $DT$-graph of $B_e'$ (and thus of $G-e$). Moreover, if $x_1 \in V(B_f)$,  then $x_1$ is a cutvertex of $G-f$ (see Theorem~G).                 

\vspace{1mm}
In either case, we reduce $G$ to the graph $H$  by replacing either $B_e'$ (in case (i)) or else $B_f$ (in case (ii)) by a path of length $3$.  By induction, $H^2$ has an ${\cal F}_4$ $x_1x_2$-hamiltonian path. This hamiltonian path can be converted to an ${\cal F}_4$ $x_1x_2$-hamiltonian path in $G^2$ by the same method used in \cite{fle2:refer}.

\vspace{2mm}
(b)  Suppose $x_1 \in V(B_e') - c_1$.

\vspace{1mm} In $(B_1)^2$, we take an $x_1c_1$-hamiltonian path $P_1(x_1,c_1)$ containing edges of $B_1$ incident to $x_r$ if $x_r$ is in $B_1$ for every $r \in \{3, 4\}$. For each $i \in \{2, \ldots, k-1\}$, we take a $c_{i-1}c_i$-hamiltonian path $P_i(c_{i-1}, c_i)$   in $(B_i)^2$ containing edges of $B_i$ incident to $x_r$ if $x_r\in V(B_i)$, for every $r \in \{3, 4\}$. In $(B_k)^2$, we take a $c_{k-1}x_2$-hamiltonian path $P_k(c_{k-1}, x_2)$. Note that this is always possible either trivially or by induction (to get an ${\cal F}_4$ hamiltonian path in $(B_i)^2$) or by Theorem \ref{strongf3} (to get a strong ${\cal F}_3$ hamiltonian path in $(B_i)^2$).Then $$P_1(x_1,c_1) \cup P_2(c_1, c_2)\cup \cdots \cup P_k(c_{k-1}, x_2)$$ yields a required $x_1x_2$-hamiltonian path in $G^2$.

\vspace{2mm} 
Consequently, if $x_1 \not \in V(B_e') - c_1$, then $x_r \in V(B_e') - c_1$ for at least one $r \in \{3, 4\}$.

\vspace{2mm}
(c) Suppose $x_1 =c_1$. 

\vspace{1mm}
Because of case (a) settled already, we have $x_r \in V(B_e') -c_1$ for at least one $r \in \{3, 4\}$.

\vspace{1mm}
(i) Suppose $x_3, x_4 \in V(B_e') - c_1$.

\vspace{1mm}
Proceeding  as in case (b), we can construct a $c_1x_2$-hamiltonian path $P_2(c_1, x_2)$ in  $(G_0 \cup B_e)^2$ containing an edge $c_1w_1$ of $B_2$.

\vspace{1mm} If there is a vertex $w \in N(c_1) \cap V(B_1)$ such that $w \not \in \{x_3, x_4\}$, then by induction let $P_1(x_1, w)$ be an ${\cal F}_4$ $x_1w$-hamiltonian path in $(B_1)^2$ containing an edge of $B_1$ incident to $x_r$ for each $r \in \{3, 4\}$. A required $x_1x_2$-hamiltonian path in $G^2$ is given by  $P_1(x_1, w) \cup (P_2(c_1, x_2) - c_1w_1) \cup \{ww_1\}$.

\vspace{1mm}
So assume that $N(c_1) \cap V(B_1) = \{x_3, x_4\}$. Then by Theorem~E let $C_1$ be an $[x_1; x_3]$-hamiltonian cycle in $(B_1)^2$ containing $x_1x_3, x_1x_4, x_3w_3$ which are edges of $B_1$. Let $P_1(x_1, x_3) = C_1 - x_1x_3$. Then a required $x_1x_2$-hamiltonian  path in $G^2$ is given by  $P_1(x_1, x_3) \cup (P_2(c_1, x_2) - c_1w_1) \cup \{x_3w_1\}$.

\vspace{1mm}
(ii) Suppose $x_3 \in V(B_e') - c_1$ and $x_4 \in V(G_0)$.

\vspace{1mm}
Let $w \in N(x_1) \cap V(B_1)$. By Theorem \ref{fletheorem4}, there is an $x_1w$-hamiltonian path $P_1(x_1, w)$ in $(B_1)^2$ containing an edge $x_3w_3$ of $B_1$.  As in case (i), we can construct an $x_1x_2$-hamiltonian path $P_2(x_1, x_2)$  in $(G_0 \cup B_e)^2$ containing an edge $x_1w_1$ of $B_2$ and an edge $x_4w_4$ of $G_0$ (apply Theorem 3 if $w_1,x_4\in V(B_2)$ and apply Theorem F otherwise). Then a required $x_1x_2$-hamiltonian path in $G^2$ is given by $P_1(x_1, w) \cup (P_2(x_1, x_2) - x_1w_1) \cup \{ww_1\}$.

\vspace{2mm}
(d) Suppose $x_1 \in V(G_0) - c_1$.

\vspace{1mm}
Then $x_1 \in V(B_t)$ for some $t \in \{2, \ldots, k-1\}$. In the case that $x_1$ is a cutvertex of $G-e$, then $x_1 = c_t$ with $t<k-1$. Let $G_t= B_1 \cup \cdots \cup B_t$ and
$H_t = B_{t+1} \cup \cdots \cup B_k$.

\vspace{2mm}
(i) Suppose $\{x_3, x_4 \} \subseteq V(G_t)$. 

\vspace{1mm}
Then by induction or by applying Theorem \ref{fletheorem4} or Theorem 3 to each 2-connected block of $G_t$, we can construct an $x_1x'$-hamiltonian path $P_1(x_1, x')$ in $(G_t)^2$ containing $x_3z_3, x_4z_4$ which are edges of $G_t$. Since $X \cap (V(H_t)-c_t) = \{x_2\}$, by applying  Theorem \ref{fletheorem3} to each 2-connected block of $H_t$, we can construct a hamiltonian cycle $C_e$ in $(H_t)^2 - c_t$ containing an edge $x_2v$ of $B_e$. 
Then a required $x_1x_2$-hamiltonian path in $G^2$ is defined by $P_1(x_1,x') \cup (C_e - x_2v) \cup \{x'v\}$.

\vspace{2mm}
(ii) Suppose $\{x_3, x_4\} \cap V(G_t) \neq \emptyset$ and $\{x_3, x_4\} \cap V(H_t) \neq \emptyset$.

\vspace{1mm}
Assume without loss of generality that $x_3 \in V(G_t)$ and $x_4 \in (V(H_t) - V(B_k))$. Because of the preceding discussion, we have $x_3 \in V(B_1)-c_1$.

\vspace{1mm}
Suppose $x_4 \in V(B_q)$ where $t < q < k$. Split $H_t$ into two block chains $J_t$ and $L_q$ where $J_t = B_{t+1} \cup \cdots \cup B_q$ and $x_4$ is not a cutvertex of $J_t$; and $L_q = B_{q+1} \cup \cdots \cup B_k$.

\vspace{1mm}
Let $P_1(x_1, x')$ denote an $x_1x'$-hamiltonian path in $(G_t)^2$ containing $x_3w_3, c_tw_t$ which are edges of $G_t$. Note that this is possible because $x_3\not= c_t$, $|V(B_t) \cap X| < 3$ and by applying Theorem 3 or Theorem F, respectively.

\vspace{1mm}
Let $C_4$ denote an hamiltonian cycle in $(J_t)^2$ containing $c_tz_t, x_4z_4$ which are edges of $J_t$. Note that this is possible by applying Theorem \ref{fletheorem3} to each block of $J_t$, provided $J_t$ is not a bridge of $H_t$. In the case that $J_t$ is a bridge $c_tx_4$, then  $C_4$ denotes  $c_tx_4 $ in $(J_t)^2$.

\vspace{1mm}
Proceed analogously to case (i) to obtain  a hamiltonian cycle $C_e$  in $(L_q)^2 - c_q$ containing an edge $x_2v$ of $B_e$.  Then a required $x_1x_2$-hamiltonian path in $G^2$ is defined by $$(P_1(x_1,x') - c_tw_t) \cup (C_4 - c_tz_t) \cup (C_e - x_2v) \cup \{w_tz_t,  x'v\}$$ if $J_t$ is not a bridge; otherwise it is defined by $$(P_1(x_1,x') - c_tw_t) \cup \{c_tx_4\} \cup (C_e - x_2v) \cup \{w_tx_4,  x'v\}.$$ This settles case (d) and thus finishes the proof of case (2.1).

\vspace{3mm}
For the remaining cases of the proof of the theorem, we adopt a different strategy of proof. For this purpose, let $B^+$ denote the graph obtained from $B_e' \cup G_0$ by adding a new edge $cx'$. Then $B^+$ is an edge-critical block. Since $|V(B^+)| < |V(G)|$, by induction, $B^+$ has the ${\cal F} _4$ property. Also, as before it is tacitly assumed that the hamiltonian paths constructed in the $(B^+)^2$ will traverse as many edges of $B^+$ as possible.

\vspace{2mm} We note that, $E((B^+)^2) = E((B_e' \cup G_0)^2) \cup E^+ $ where
\[   E^+ = \{cx'\} \cup \{  x'w_c, u'c \ | \ w_c \in N(c) \cap V(G_0), u' \in N(x') \cap V(B_e')\}.     \]

In what follows, any vertex in $N(c) \cap V(G_0)$ will be subscribed with $c$, and any vertex in $N(x') \cap V(B_e')$ will be superscribed with $'$. Also, we use  $y$  to denote a neighbor of $x$ in $B_e$.

\vspace{2mm}  (2.2) Suppose $x_2 \in V(B_e)$ and $x_2 \neq x$.

\vspace{2mm} Then $x_1, x_3, x_4 \in V(B^+)-\{c\}$. Let $P^+(x_1, c)$ denote an ${\cal F} _4$ $x_1c$-hamiltonian path in $(B^+)^2$ containing $x_3z_3, x_4z_4$ which are different edges of $B^+$ using induction. Note that $x_iz_i=x'c$ is possible for $i \in\{3,4\}$.

\vspace{2mm}
Set $E^*=E(P^+(x_1, c)) \cap E^+$ and set $|E^*| = r$. Clearly, $0 \leq r \leq 3$. Observe that $r=4$ would imply that $x'$ and $c$ are internal vertices of the corresponding hamiltonian path, which is not possible. However, the case $r=3$ could be reduced to the case (b) (i) below traversing more edges of $B^+$ than the original path. Thus $r=3$ is also impossible.

\vspace{2mm}
(a) Suppose  $r=0$ in which case $x_iz_i\neq x'c$ for $i=3,4$.

\vspace{1mm}
Trivially $(B_e)^2$ has an $x_2c$-hamiltonian path $P_2(x_2, c)$ (since $x_2 \notin \{c, x\}$). Then  $P^+(x_1, c) \cup P_2(x_2, c)$ defines a required  ${\cal F} _4$ $x_1x_2$-hamiltonian path in $G^2$.

\vspace{2mm} 
(b) Suppose   $r=1$.    

\vspace{1mm}
By Theorem \ref{fletheorem4} {\it (i)}, $(B_e)^2$ has an $x_2x$-hamiltonian path $P_2(x_2, x)$ containing an edge $cw$ of $B_e$.

\vspace{1mm}
(i) If $E^* = \{cx'\}$, then $(P^+(x_1, c) -  cx')  \cup P_2(x_2, x)  \cup \{xx'\}$ defines a required  ${\cal F} _4$ $x_1x_2$-hamiltonian path in $G^2$. Note  that $x_iz_i=x'c$ for $i\in \{3,4\}$ is not an obstacle.

\vspace{1mm}
From now on we ca assume that $x_iz_i\neq x'c$ for $i=3,4$.

\vspace{1mm}
(ii) If $E^*= \{cu'\}$, then
$(P^+(x_1, c) -  cu')  \cup P_2(x_2, x)  \cup \{xu'\}$ defines a required  ${\cal F} _4$ $x_1x_2$-hamiltonian path in $G^2$. 

\vspace{1mm}
(iii) If $E^* = \{x'w_c\}$, then
$(P^+(x_1, c) -  x'w_c)  \cup (P_2(x_2,x) - cw)  \cup \{xx', ww_c\}$ defines a required  ${\cal F} _4$ $x_1x_2$-hamiltonian path in $G^2$. This holds true even if $w_c\in \{x_3,x_4\}$ and $w_cc\in E(P^+(x_1,c))$.

\vspace{2mm} (c) Suppose   $r=2$.

\vspace{1mm}
By Theorem \ref{fletheorem4},  $(B_e)^2$ has an $x_2c$-hamiltonian path $P_2(x_2, c)$ containing an edge $xy$ of $B_e$.

\vspace{1mm}
(i) If $E^* = \{cx', x'w_c\}$, then
$$ (P^+(x_1, c) - \{cx', x'w_c\})  \cup (P_2(x_2, c) - xy) \cup \{yx', xx', cw_c\} $$  defines a required  ${\cal F} _4$ $x_1x_2$-hamiltonian path in $G^2$, even if $\{x_3,x_4\}\cap\{x',w_c\}\not=\emptyset$. Note  that $x_iz_i=x'c$ for $i\in \{3,4\}$ is not an obstacle.

\vspace{1mm}
From now on we ca assume that $x_iz_i\neq x'c$ for $i=3,4$.

\vspace{1mm}
(ii)  If $E^* = \{x'u_c, x'w_c\}$, then
$$(P^+(x_1, c) - \{x'u_c, x'w_c\})  \cup (P_2(x_2, c) - xy) \cup \{ yx', xx', w_cu_c\} $$ defines a required  ${\cal F} _4$ $x_1x_2$-hamiltonian path in $G^2$.

\vspace{2mm}
(iii) If $E^* = \{x'u_c, w'c\}$, then
$$  (P^+(x_1, c) - \{x'u_c, w'c\})  \cup (P_2(x_2, c)) \cup \{cu_c, w'x'\}  $$  defines a required  ${\cal F} _4$ $x_1x_2$-hamiltonian path in $G^2$.


\vspace{3mm}  (2.3) Suppose $x_4 \in V(B_e)$.

\vspace{2mm} Then $x_1, x_2, x_3 \in V(B^+)-\{c\}$. Let $P^+(x_1, x_2)$ denote an ${\cal F} _4$ $x_1x_2$-hamiltonian path in $(B^+)^2$ containing $x_3z_3, cc^*$, where $c^*\in \{w_c,x'\}$, which are different edges of $B^+$ using induction. Note that $x_3z_3=x'c$ is possible.

\vspace{2mm}
Now set $E^*=E(P^+(x_1, x_2)) \cap E^+$ and set $|E^*| =r$. Clearly, $0\leq r\leq 4$. Note that the case $r=4$ yields a contradiction just as did the case $r=3$ in the subcase (2.2) above.

\vspace{2mm}
(a)  Suppose $r=0$, in which case $c^*=w_c$ and $x_3z_3\neq x'c$.

\vspace{2mm}
By Theorem \ref{fletheorem3}, $(B_e)^2$ has a $[c; x_4]$-hamiltonian cycle $C_e$ containing $cw, cz, x_4z_4$ which are edges of $B_e$. Then $(P^+(x_1, x_2) - cw_c) \cup (C_e - cw) \cup \{ww_c\}$ defines a required  ${\cal F} _4$ $x_1x_2$-hamiltonian path in $G^2$.

\vspace{2mm}
(b) Suppose $r=1$.

\vspace{2mm}
(i) Suppose $E^* = \{x'c \}$ or $E^* = \{w'c \}$. By Theorem \ref{fletheorem4}, $(B_e)^2$ has a $cx$-hamiltonian path $P_4(c,x)$ containing an edge $x_4z_4$ of $B_e$. Then $(P^+(x_1, x_2) - zc) \cup P_4(c,x) \cup \{zx\}$ is a required  ${\cal F} _4$ $x_1x_2$-hamiltonian path in $G^2$ for any vertex  $z \in \{x', w'\}$. Note that either $c^*=x'$ or $x_3z_3=x'c$ is not an obstacle. 

\vspace{2mm}
(ii) Suppose $E^* = \{x'u_c \}$. Hence $c^*=w_c$ and $x_3z_3\neq x'c$. By Theorem \ref{strongf3}, $(B_e)^2$ has a $cx$-hamiltonian path $P_4(c,x)$ containing $cw, x_4z_4$ which are edges of $B_e$. As for the $x_1x_2$-hamiltonian path $P^+(x_1, x_2)$ in $(B^+)^2$ we possibly have $u_c =w_c$. In any case,
$(P^+(x_1, x_2) - x'u_c) \cup (P_4(c,x) - cw) \cup \{x'x, wu_c\}$ defines a  required  ${\cal F} _4$ $x_1x_2$-hamiltonian path in $G^2$.

\vspace{2mm}
(c) Suppose $r=2$. Note that $x_4\neq c$ in this case and $xc\notin E(B_e)$ because of $G$ is edge-critical.

\vspace{2mm}
(c1) $E^* = \{x'c, w'c \}$, in which case $c^*=x'$ and hence $x_3z_3\neq x'c$. By Theorem \ref{fletheorem4}, $(B_e)^2$ has an $xy$-hamiltonian path $P_4(x,y)$ containing an edge $x_4z_4$ of $B_e$. Then $(P^+(x_1, x_2) - \{x'c, w'c \}) \cup P_4(x, y) \cup \{x'y, w'x\}$ yields a required  ${\cal F} _4$ $x_1x_2$-hamiltonian path in $G^2$.

\vspace{2mm}
(c2) $E^* = \{x'c, x'u_c \}$. 

\vspace{2mm}
Let $y \in N(x) \cap V(B_e)$ where $y \neq x_4$, and let $P_4(x, y)$ be an ${\cal F} _4$ $xy$-hamiltonian path in $(B_e)^2$ containing  $x_4z_4, cw$ which are edges  of $B_e$ (by induction or by Theorem 3). Then $(P^+(x_1, x_2) - \{x'c, x'u_c \}) \cup (P_4(x, y) - cw) \cup \{x'x, x'y, wu_c\}$  results in a required  ${\cal F} _4$ $x_1x_2$-hamiltonian path in $G^2$. Note that either $c^*=x'$ or $x_3z_3=x'c$ is not an obstacle.

\vspace{2mm}
From now on we can assume that $c^*=w_c$ and $x_3z_3\neq x'c$.

\vspace{2mm}
(c3) $E^* = \{x'y_c, x'u_c \}$. 

\vspace{2mm}  
By Lemma 3, there is a $[c; x, x_4]$-hamiltonian cycle $C_e$ in $(B_e)^2$ containing $cw$, $cu$, $xy$, $x_4z_4$ which are edges of $B_e$ provided $x_4\neq x$; otherwise, let $C_4$ be an $[c;x_4]$-hamiltonian cycle of $(B_e)^2$ containing $cw, cu, xy$ which are edges of $B_e$ resulting from an application of Theorem E. Then $$(P^+(x_1, x_2)-\{x'y_c, x'u_c \})\cup (C_e-\{cw, cu, xy\}) \cup \{wy_c, uu_c, xx', yx'\}$$ defines a required  ${\cal F} _4$ $x_1x_2$-hamiltonian path in $G^2$ independent of the position of $x_4$. 


\vspace{2mm}
(c4) $E^* = \{x'y_c, w'c \}$. Then there are two subcases to consider.

\vspace{2mm}
(i) Suppose $y_c = w_c$.   Let $P_4(x, y)$ be as defined in case (c2). Then $$(P^+(x_1, x_2) - \{x'y_c, w'c \}) \cup (P_4(x, y) - cw) \cup \{w'x, x'y, wy_c \}$$  yields a required  ${\cal F} _4$ $x_1x_2$-hamiltonian path in $G^2$.

\vspace{2mm}
(ii) Suppose $y_c \neq w_c$. There are four possibilities.

\vspace{1mm} If $P^+(x_1, x_2)$ takes the form $x_1 \cdots x'y_c \cdots w_ccw' \cdots x_2$, then proceed as in (i) to obtain  a required  ${\cal F} _4$ $x_1x_2$-hamiltonian path in $G^2$.

\vspace{1mm} If $P^+(x_1, x_2)$ takes the form $x_1 \cdots y_cx' \cdots  w_ccw' \cdots x_2$\ or $x_1 \cdots x'y_c \cdots  w'cw_c \cdots x_2$, then we can reduce this case to case (a) where $r=0$ as follows. Delete $x'y_c, w'c$ from $P^+(x_1, x_2)$ and add to it the edges $x'w', cy_c$.

\vspace{1mm} If $P^+(x_1, x_2)$ takes the form $x_1 \cdots y_cx' \cdots  w'cw_c \cdots x_2$, then let $P_4(x, y)$ denote an $xy$-hamiltonian path in $(B_e)^2$ as defined in case (c2). Then
$$(P^+(x_1, x_2) - \{x'y_c, w'c \}) \cup (P_4(x, y) - cw) \cup \{w'x, x'y, y_cw\}$$  defines a required
${\cal F} _4$ $x_1x_2$-hamiltonian path in $G^2$.

The other cases are symmetrical. 

\vspace{2mm}
(c5) $E^*=\{cv', cw'\}$. This case cannot happen since $E^*\cap E(B^+)=\emptyset$, but $cw_c\in E(P_2(x_1,x_2)) \cap E(B^+)$.

(d) Suppose $r=3$. Thus $E^*$ must be one of the following three sets: $E^*=\{x'u_c, x'c, w'c\}$, $E^*=\{x'u_c, x'y_c, w'c\}$,
$E^*=\{w'c, v'c, x'u_c\}$. It is now straightforward to see that in each of these three cases the corresponding $P^+(x_1,x_2)$ can be modified so as to contain more edges of $B^+$  and satisfying $E^*=\{x'c\}$, i.e., $r=1$. Namely, in the respective case 

form $(\{x'u_c,x'c,w'c\}-\{x'u_c,w'c\})\cup \{w'x',cu_c\}$;

replace $\{x'u_c,x'y_c,w'c\}$ with $\{w'x',x'c,u_cy_c\}$;

replace $\{w'c,v'c,x'u_c\}$ with $\{w'v', x'c, cu_c\}$.

\vspace{2mm}   Theorem 4 now follows.                                              \qed

 \vspace{2mm} As a special case of Theorem 4 we obtain the following.
 
 \begin{deduce} 
  \vspace{1mm} Let $G$ be a $2$-connected graph on $n\geq 4$ vertices, and let $e=xy\in E(G)$ and $u,v\in V(G)$ such that $\{x,y\}\cap\{u,v\}=\emptyset$. Then $G^2$ has a hamiltonian cycle $C$ with $e\in E(C)$, and at least one of the edges of $u$ in $C$ at least one of the edges of $v$ in $C$ are edges of $G$.
  \end{deduce}

\vspace{2mm}

\section{Final remarks}
In subsequent papers we shall use some of the theorems of this paper to describe (among other results) the most general structure a graph may have such that its square is hamiltonian or hamiltonian connected, respectivelly. This will also solve a problem raised in \cite{cot:refer} in the affirmative and proves a conjecture raised in \cite{tra:refer}; we shall also present a partial solution of a conjecture stated in \cite{eks:refer}.

\vspace{1mm}
 It is easy to see that the complete bipartite graph $K_{2,k-2}$  does not have the ${\cal F} _k$ property for every integer $k \geq 5$. For example, take $x_1, x_2$ to be the two vertices of degree $k-2$ and $x_3, \ldots, x_k$ to be the rest of the vertices.  Hence Theorem \ref{f4} is best possible.

\vspace{2mm}
A graph $G$ is said to have the {\em ${\cal \overline{F}}$ property} if it has three $2$-valent  vertices $x, y, z$ such that $N(x) =N(y) = N(z)$. From the above observation, we see that if $G$ has the ${\cal \overline{F}} $ property, then $G$ does not have the ${\cal F} _k$ property for any $k \geq 5$.

\vspace{2mm}
While it is now known that Theorem \ref{fletheorem4}(i) can be generalized to Theorem \ref{f4}, it is also of interest to know whether or not Theorem \ref{fletheorem3} can be generalized to $3$ given vertices. That is, given three arbitrary vertices $v, w_1, w_2$ of a $2$-connected graph $G$, does $G^2$ contain a $[v;w_1,w_2]$-hamiltonian cycle $C$? The following example shows that this is not true in general.

\vspace{2mm}
Let $k\geq 5$ be an integer  and let $v_1v_2\cdots v_nv_1$ be a  cycle with $n$ vertices where $n \geq k+3$. Take a new vertex $v$ and join it to $v_1$ and $v_k$ to get the graph $H$. Let $w_1 = v_1$ and $w_2=v_k$. Then it is easy to see that $H^2$ admits  no hamiltonian cycle $C$ containing the edges $vw_1, vw_2$ and $w_iz_i$ where $z_i \in N(w_i)$, $i=1,2$.

\vspace{5mm}
\begin{center}    {\bf Acknowledgements} \end{center}

Research of the first author was supported in part by FWF-grant P27615-N25, whereas research of the second  author was supported by the $FRGS$ Grant (FP036-2013B).

\end{document}